\documentclass[11pt]{amsart}
\input epsf
\usepackage{latexsym}
\usepackage{amssymb,amsmath,amsthm,amscd, amssymb,
graphicx, enumerate, verbatim}
\usepackage[all]{xy}

\numberwithin{equation}{section}
\newtheorem{cor}[equation]{Corollary}
\newtheorem{lem}[equation]{Lemma}
\newtheorem{prop}[equation]{Proposition}
\newtheorem{thm}[equation]{Theorem}

\newtheorem{defn}[equation]{Definition}

\newtheorem{Example}[equation]{Example}

\newtheorem{remark}[equation]{Remark}
\newenvironment{rmk}{\begin{remark}\rm}{\end{remark}}
\def\co{\colon\thinspace}
\def\conv#1{{\underset{#1}{\longrightarrow}}}

\newcommand{\curv}{\mbox{curv}}

\newcommand{\vol}{\mbox{vol}}
\newcommand{\diam}{\mbox{diam}}

\newcommand{\Iso}{\mbox{Iso}}

\newcommand{\e}{\epsilon}
\def\a{\alpha}
\def\G{\Gamma}

\def\d{\delta}

\def\s{\sigma}
\def\l{\lambda}
\def\L{\Lambda}
\def\Z{\mathbb Z}
\def\Rn{\mathbb{R}^n}
\def\Rk{\mathbb{R}^k}
\newcommand{\g}{\gamma}

\begin{document}

\abovedisplayskip=6pt plus3pt minus3pt
\belowdisplayskip=6pt plus3pt minus3pt

\title[Negatively pinched manifolds with amenable fundamental groups ]{\bf
Classification of negatively pinched\\ manifolds with amenable
fundamental groups}

\rm

\author{Igor Belegradek}
\address{Igor Belegradek\\School of Mathematics\\ Georgia Institute of
Technology\\ Atlanta, GA 30332-0160}\email{ib@math.gatech.edu}
\author{Vitali Kapovitch}
\address{Vitali Kapovitch\\
Department of Mathematics\\
University of Maryland\\
College Park, MD 20742, USA}
\email{vtk@math.umd.edu}
%
\date{}
\begin{abstract}
We give a diffeomorphism classification of pinched negatively
curved manifolds with amenable fundamental groups, namely, they
are precisely the M\"obius band, and
the products of $\mathbb R$ with  the total spaces
of flat vector bundles over closed infranilmanifolds.
\end{abstract}
\maketitle

\section{Introduction}

In this paper we study  manifolds of the form $X/\G$, where $X$ is
a simply-connected complete Riemannian manifold with sectional
curvatures pinched (i.e. bounded) between two negative constants,
and $\G$ is a discrete torsion free subgroup of the isometry group
of $X$. According to~\cite{BuSc}, if
$\G$ is amenable, then either $\G$ stabilizes a biinfinite
geodesic, or else $\G$ fixes a unique point $z$ at infinity. The
case when $\G$ stabilizes a biinfinite geodesic is completely
understood, namely the normal exponential map to the geodesic is a
$\G$-equivariant diffeomorphism, hence $X/\G$ is a vector bundle
over $S^1$; there are  only two such bundles each admitting a
complete hyperbolic metric.

If $\G$ fixes a unique point $z$ at infinity (such groups are
called {\it parabolic}), then $\G$ stabilizes horospheres centered at
$z$ and permutes geodesics asymptotic to $z$, so that given a
horosphere $H$, the manifold $X/\G$ is diffeomorphic to the
product of $H/\G$ with $\mathbb R$. We refer to $H/\G$ as a {\it
horosphere quotient}. In this case a delicate result of
B.~Bowditch~\cite{Bow} shows that
$\G$ must be finitely generated, which by Margulis
lemma~\cite{BGS} implies that $\G$ is virtually nilpotent.

The main result of this paper is a diffeomorphism classification
of horosphere quotients, namely we show that, up to a
diffeomorphism, the classes of horosphere quotients and (possibly
noncompact) infranilmanifolds coincide.

By an {\it infranilmanifold} we mean the quotient of a
simply-connected nilpotent Lie group $G$ by the action of a
torsion free discrete subgroup $\G$ of the semidirect
product of $G$ with a compact subgroup of $\mathrm{Aut}(G)$.
\begin{thm}\label{main thm}
For a smooth manifold $N$ the following are equivalent:\newline
$\mathrm{\ \ \ \ (1)}$ $N$ is a horosphere quotient;\newline
$\mathrm{\ \ \ \ (2)}$ $N$ is diffeomorphic to an infranilmanifold;\newline
$\mathrm{\ \ \ \ (3)}$ $N$ is the total space of a flat
Euclidean vector bundle over a compact infranilmanifold.
\end{thm}
The implication $(3)\Rightarrow (2)$ is straightforward,
$(2)\Rightarrow (1)$ is proved by constructing an explicit warped
product metric of pinched negative curvature. The proof of
$(1)\Rightarrow (3)$ occupies most of the paper, and depends on
the collapsing theory of J.~Cheeger, K.~Fukaya, and
M.~Gromov~\cite{CFG}.

If $N$ is compact (in which case the conditions
$\mathrm{(2)}$, $\mathrm{(3)}$ are identical),
the implication $(1)\Rightarrow (2)$
follows from Gromov's classification of almost flat
manifolds, as improved by E.~Ruh, while the implication
$(2)\Rightarrow (1)$ is new.
If $N$ is non-compact, Theorem~\ref{main thm}
is nontrivial even when $\pi_1(N)\cong\mathbb Z$, although
the proof does simplify in this case.
A direct algebraic proof of $(2)\Rightarrow (3)$ was given
in~\cite[Theorem 6]{Wil2}, but the case when
$N$ is a nilmanifold was already treated in~\cite{Mal}, where
it is shown that any nilmanifold is diffeomorphic to
the product of a compact nilmanifold and a Euclidean space.

We postpone the discussion of
the proof till Section~\ref{sec: ideas}, and just mention that
the proof also gives geometric information about horosphere
quotients, e.g. we show that $H/\G$ is diffeomorphic to a tubular
neighborhood of some orbit of an N-structure on $H/\G$.

By Chern-Weil theory 
any flat Euclidean vector bundle has zero rational Euler and
Pontrjagin classes. Moreover,
by~\cite{Wil1} any flat Euclidean bundle with virtually abelian
holonomy is isomorphic to a bundle with finite structure group. 
Thus the vector bundle in $(3)$ becomes trivial in a finite cover, and
has zero rational Euler and Pontrjagin classes,
and in particular, any horosphere quotient is finitely covered
by the product of a compact nilmanifold and a Euclidean space.

\begin{cor}\label{intro: maincor}
A smooth manifold $M$ with amenable fundamental group admits a
complete metric of pinched negative curvature if and only if it is
diffeomorphic to the M\"obius band, or to the product of a line
and the total space a flat Euclidean vector bundle over a compact
infranilmanifold.
\end{cor}

The pinched negative curvature assumption in Corollary~\ref{intro:
maincor} cannot be relaxed to $-1\le \sec\le 0$ or $\sec\le -1$,
e.g. because these assumptions do not force the fundamental group
to be virtually nilpotent~\cite[Section 6]{Bow}. More delicate
examples come from the work of M.~Anderson~\cite{And} who proved
that each vector bundle over a closed nonpositively curved
manifold (e.g. a torus) carries a complete Riemannian metric with
$-1\le \sec\le 0$. Since in each dimension there are only finitely
many isomorphism classes of flat Euclidean bundles over a given
compact manifold, all but finitely many vector bundles over tori
admit no metrics of pinched negative curvature. Also $-1\le
\sec(M)\le 0$ can be turned into $\sec(M\times\mathbb R)\le -1$
for the warped product metric on $M\times\mathbb R$ with warping
function $e^t$~\cite{BON}, hence Anderson's examples carry metrics
with $\sec\le -1$ after taking product with $\mathbb R$.
Specifically, if $E$ is the total space of a vector bundle over a
torus with nontrivial rational Pontrjagin class, then $M=E\times
\mathbb R$ carries a complete metric of $\sec\le -1$ but not of
pinched negative curvature. Finally, Anderson also showed that
every vector bundle over a closed negatively curved manifold
admits a complete Riemannian metric of pinched negative curvature,
hence amenability of the fundamental group is indispensable.

Because an infranilmanifold with virtually abelian
fundamental group is flat, Theorem~\ref{main thm} immediately implies
the following.
\begin{cor}\label{cor:ab}
Let $M$ be a smooth manifold with virtually abelian fundamental group.
Then the following are equivalent\newline
$\mathrm{\ \ \ \ (i)}$ $M$ admits a complete metric of $\sec\equiv -1$;\newline
$\mathrm{\ \ \ \ (ii)}$ $M$ admits
a complete metric of pinched negative curvature.
\end{cor}

In~\cite{Bow2} Bowditch developed several equivalent definitions
of geometrical finiteness for pinched negatively curved manifolds,
and conjectured the following result.

\begin{cor}\label{cor: geom finite}
Any geometrically finite pinched negatively curved manifold $X/\G$
is diffeomorphic to the interior of a compact manifold with
boundary.
\end{cor}

We believe that the main results of this paper, including
Corollary~\ref{cor: geom finite},  should extend to the orbifold
case, i.e. when $\G$ is not assumed to be torsion
free. However, working in the orbifold category
creates various technical
difficulties, both mathematical and expository, and we do not
attempt to treat the orbifold case in this paper.

\bf{Acknowledgements:}\rm\ It is a pleasure to thank Ilya Kapovich
and Derek J. S. Robinson for Lemma~\ref{root}, Hermann Karcher for
pointing us to~\cite{BrKa}, Robion C. Kirby for
Lemma~\ref{l:isotopy}, Anton Petrunin for  numerous helpful
conversations and suggestions and  particularly for
Lemma~\ref{crit}, and Xiaochun Rong and Kenji Fukaya for helpful
discussions on collapsing. We are grateful to the referee
for various editorial comments.
This work was partially supported by the NSF grants
\# DMS-0352576 (Belegradek) and \# DMS-0204187 (Kapovitch).

\tableofcontents

\section{Sketch of the proof of $(1)\Rightarrow (3)$}
\label{sec: ideas}

The Busemann function corresponding to $z$ gives rise to a
$C^2$-Riemannian submersion $X/\G\to \mathbb R$ whose fibers
are horosphere quotients each equipped with the induced
$C^1$-Riemannian metric $g_t$. By the Rauch comparison theorem
the second fundamental form of a horosphere is bounded in terms
of curvature bounds of $X$ (cf.~\cite{BrKa}).
In particular, each fiber has curvature uniformly bounded above
and below in comparison sense. Let $\s(t)$ be a horizontal
geodesic in $X/\G$, i.e. a geodesic that projects isometrically
to $\mathbb R$.
Because of the exponential convergence of geodesics
in $X$, the manifold $X/\G$ is ``collapsing'' in the sense that
the unit balls around $\s(t)\in X/\G$ form an
exhaustion of $X/\G$ and have small injectivity radius for large $t$.
Similarly, each fiber of $X/\G\to \mathbb R$ also collapses, and in fact
$X/\G$ is noncollapsed in the direction transverse to the fibers.

There are essential difficulties in applying the
collapsing theory of~\cite{CFG} to $X/\G$. First, we do not know
whether $(X/\G, \s(t))$ converges in pointed Gromov-Hausdorff
topology to a single limit space. By general theory, the family
$(X/\G, \s(t))$ is precompact and thus has many converging
subsequences. While different limits might be
non-isometric, one of the main steps of the proof is obtaining a uniform
(i.e. independent of the subsequence) lower bound on the
``injectivity radius'' of the limit spaces at the base point. This
is done by a comparison argument involving taking ``almost square
roots'' of elements of $\G$, and using the flat connection
of~\cite{Bow} discussed below.
Another complication is that
the N-structure on $(X/\G, \s(t))$ provided by~\cite{CFG}
may well have zero-dimensional orbits outside the unit ball around
$\s(t)$, in other words a large noncompact region of $(X/\G, \s(t))$
may be non-collapsed, which makes it hard to control topology of
the region.

However, once the ``injectivity radius'' bound is established,
critical points of distance functions considerations
yield the ``product structure at infinity'' for $X/\G$, and also
for $(H/\G, g_t)$ if $t$ is large enough.

Furthermore, one can show that $H/\G$ is diffeomorphic to the
normal bundle of an orbit $O_t$ of the N-structure. The orbit
corresponds to the point in a limit space given to us by the
convergence, and at which we get an ``injectivity radius'' bound.
This depends on a few results on Alexandrov spaces with curvature
bounded below, with  key ingredients provided by~\cite{Kap}
and~\cite{PP}.

By the collapsing theory, the structure group of the normal bundle
to the orbit of an N-structure is a finite extension of a
torus group~\cite{CFG}.
Of course, not every such a bundle has a flat Euclidean structure.

The flatness of the normal bundle to the orbit is proved using a
remarkable flat connection discovered by B.~Bowditch~\cite{Bow},
and later in a different disguise by W.~Ballmann and
J.~Br\"uning~\cite{BB}, who were apparently unaware of~\cite{Bow}.
It follows from~\cite{Bow} or~\cite{BB} that each pinched
negatively curved manifold $X/\G$, where $\G$ fixes a unique point
$z$ at infinity, admits a natural flat $C^0$-connection that is
compatible with the metric and has nonzero torsion, and such that
on short loops it is close to the Levi-Civita connection.
Furthermore, the parallel transport of the connection preserves
the fibration of $X/\G$ by horosphere quotients. Hence each
horosphere quotient has flat tangent bundle.

In fact, we prove a finer result
that the normal bundle to $O_t$ is also flat, for suitable large $t$.
(Purely topological considerations are useless here since there
exist vector bundles without flat Euclidean structure
whose total spaces have flat Euclidean tangent bundles, for example this
happens for any nontrivial orientable $\mathbb R^2$-bundle over
the $2$-torus that has even Euler number).
It turns out that $O_t$ sits with flat normal
bundle in a totally geodesic stratum of the N-structure, so it
suffices to show that the normal bundle to the stratum is flat,
when restricted to $O_t$. Now since the above flat connection is
close to Levi-Civita connection the normal bundle is ``almost
flat'', and it can be made flat by averaging via center of mass.
This completes the proof.

Throughout the proof we use
the collapsing theory developed in~\cite{CFG}.
This paper is based on the earlier extensive work of Fukaya,
and Cheeger-Gromov, and many arguments in~\cite{CFG} are merely sketched.
We suggest reading~\cite{Fukaya-survey} for a snapshot of
the state of affairs before~\cite{CFG}, and~\cite{PT, PRT, FR}
for a current point of view.

\section{Topological digression}
The result of Bowditch~\cite{Bow} that horosphere quotients
have finitely generated fundamental groups actually
implies that
any horosphere quotient is homotopy equivalent to a compact
infranilmanifold (because any torsion free finitely generated
virtually nilpotent group is the fundamental group of a compact
infranilmanifold~\cite{Dek}, and because for aspherical manifolds
any $\pi_1$-isomorphism is induced by a homotopy equivalence).

To help appreciate the difference between this statement
and Theorem~\ref{main thm}, we discuss several types of
examples that are allowed by Bowditch's result, and are
ruled out by Theorem~\ref{main thm}. The simplest example
is a vector bundle over an infranilmanifold with nonzero
rational Euler or Pontrjagin class: such a manifold cannot be
the total space of a flat Euclidean bundle, as is easy to see
using the fact that the tangent bundle to any infranilmanifold
is flat.

Another example is the product of a closed infranilmanifold
and a contractible manifold of dimension $>2$
that is not simply-connected at infinity.
Finally, even more sophisticated examples come from the fact that
below metastable range (starting at which any homotopy equivalence
is homotopic to a smooth embedding, by Haefliger's embedding
theorem) there are many smooth manifolds that are thickenings of
say a torus, yet are not vector bundles over the torus. It would
be interesting to see whether the ``weird'' topological
constructions of this paragraph can be realized geometrically,
even as nonpositively curved manifolds.

\section{Parallel transport through infinity and
rotation homomorphism}
\label{sec: rot hom}

Let $X$ be a simply-connected complete pinched negatively curved
$n$-manifold normalized so that $-a^2\le sec(X)\le -1$. One of the
key properties of $X$ used in this section is that any two
geodesic rays in $X$ that are asymptotic to the same point at
infinity converge exponentially, i.e.  for any asymptotic rays
$\g_1(t), \g_2(t)$ with $\g_1(0), \g_2(0)$ lying on the same
horosphere, the function $d(\g_1(t),\g_2(t))$ is monotonically
decreasing as $t\to\infty$ and
\begin{equation*}\label{asym}
e^{-at} c_1(a, d(\g_1(0),\g_2(0)))\le d(\g_1(t),\g_2(t))\le
e^{-t} c_2(a, d(\g_1(0),\g_2(0))).
\end{equation*}
where $c_i(a,d)$ is linear in $d$ for small $d$. This is proved by
triangle comparison with spaces of constant negative curvature.

Bowditch introduced a connection on $X$ that we now
describe (see~\cite[Section 3]{Bow} for details). Fix a point $z$
at infinity of $X$. Let $w_i\to z$ as $i\to\infty$. For any $x,
y\in X$, consider the parallel transport map from $x$ to $w_i$
followed by the parallel transport from $w_i$ to $y$ along the
shortest geodesics. This defines an isometry between the tangent
spaces at $x$ and $y$. By~\cite[Lemma 3.1]{Bow}, this map
converges to a well-defined limit isometry $P^\infty_{xy}\co
T_{x}M\to T_yM$ as $i\to\infty$. We refer to $P^\infty_{xy}$ as
the {\it parallel transport through infinity from $x$ to $y$}.

We denote the Levi-Civita parallel transport from $x$ to $y$ along
the shortest geodesic by $P_{xy}$; clearly, if $x, y$ lie on a
geodesic ray that ends at $z$, then $P^\infty_{xy}=P_{xy}$. A key
feature of $P^\infty$ is that it approximates the Levi-Civita
parallel transport on short geodesic segments (see~\cite[Lemma
3.2]{Bow}; more details can be found in~\cite[Section
6]{BuKa}). This is because any geodesic triangle in $X$
spans a ``ruled'' surface of area at most the area of the
comparison triangle in the hyperbolic plane of $\sec=-1$. By
exponential convergence of geodesics the area of the comparison
triangle is bounded above by a constant times the shortest side of
the triangle. As the holonomy around the circumference of the
triangle is bounded by the integral of the curvature over its
interior, we conclude that $|P_{xy}-P^\infty_{xy}|\le q(a)d(xy)$,
where $q(a)$ is a constant depending only on $a$.

Given $x\in X$, fix an isometry $\Rn\to T_x X$ and translate it
around $X$ using $P^\infty$. This defines a $P^\infty$-invariant
trivialization of the tangent bundle to $X$. Let $\Iso_z(X)$ be
the group of isometries of $X$ that fixes $z$. For any point $y\in
X$ look at the map $\Iso_z(X)\to O(n)$ given by $\gamma\mapsto
P^{\infty}_{\g(y)y}\circ d\gamma$. It turns out that this map is a
homomorphism independent of $y$. We call it the {\it rotation
homomorphism}. Starting with a different base point $x\in X$ or a
different isometry $\Rn\to T_x X$ has the effect of replacing the
rotation homomorphism by its conjugate.

Now if $\G$ is a discrete torsion free subgroup of $\Iso_z(X)$,
then, since the rotation homomorphism is independent of $y$,
$P^\infty$ gives rise to a flat connection on $X/\G$ with holonomy
given by the rotation homomorphism. By the above discussion of
$P^\infty$, this is a $C^0$ flat connection that is compatible
with the metric and close to the Levi-Civita connection on short
loops. Of course, this connection has torsion.

\begin{rmk}\label{rmk: valid for c1}
The above discussion is easily seen to be  valid if $X$ is a
simply-connected complete $C^1$-Riemannian manifold of pinched
negative curvature in the comparison sense. This is because any
such $C^1$-metric can be approximated uniformly in $C^1$-topology
by smooth Riemannian metrics of pinched negative
curvature~\cite{Nik}, perhaps with slightly larger pinching. Then
the distance functions and Levi-Civita connections converge
uniformly in $C^0$-topology, and we recover all the statements
above.
\end{rmk}
\begin{rmk}
The connection of Bowditch, that was described above, was
reinvented later in a different disguise by Ballmann and
Br\"uning~\cite[Section 3]{BB}.
The connection in~\cite{BB} is defined by an explicit
local formula in terms of the curvature tensor and the Levi-Civita
connection of $X$. Actually,~\cite{BB} only discusses the case of
compact horosphere quotients however all the arguments there are
local, hence they apply to any horosphere quotient. The only
feature  which is special for compact horosphere quotients is that
in that case the connection has finite holonomy group~\cite{BB},
as follows from estimates in~\cite{BuKa}. For non-compact horosphere
quotients, the holonomy
need not be finite as seen by looking at a glide rotation with
irrational angle in $\mathbb R^3$, thought of as a horosphere in
the hyperbolic $4$-space. We never have to use~\cite{BB} in this
paper, however, for completeness we discuss their construction in
Appendix~\ref{BBB}, where we also show that the connections
of~\cite{Bow} and~\cite{BB} coincide.
\end{rmk}

\section{Passing to the limit}
\label{sec: pass to lim}

Let $X$ be a simply-connected complete pinched negatively curved
$n$-manifold normalized so that $-a^2\le sec(X)\le -1$, let $c(t)$
be a biinfinite geodesic in $X$, and let $\G$ be a closed subgroup
of $\Iso(X)$ that fixes the point $c(\infty)$ at infinity. We
refer to the gradient flow $b_t$ of a Busemann function for $c(t)$
as {\it Busemann flow}. Following Bowditch, we sometimes use the
notation  $x+t:=b_t(x)$.

Since
$X$ has bounded curvature and infinite injectivity radius, the
family $(X, c(t),\G)$ has a subsequence $(X, c(t_i),\G)$ that
converges to $(X_\infty, p, G)$ in the equivariant pointed
$C^{1,\a}$-topology~\cite[Chapter 10]{Pet}. 
Here $X_\infty$ is a smooth manifold with
$C^{1,\a}$-Riemannian metric that has infinite injectivity radius
and the same curvature bounds as $X$ in the comparison sense, and
$G$ is a closed subgroup of $\Iso(X_\infty)$. Note that
$\Iso(X_\infty)$ is a Lie group that acts on $X$ by
$C^3$-diffeomorphisms (this last fact  is probably known  but for
a lack of reference  we give a simple proof in
Appendix~\ref{apendix: iso are c3}).

Furthermore, geodesic rays in $X$ that start at uniformly bounded
distance from $c(t_i)$ converge to rays in $X_\infty$. In
particular, the rays $c(t+t_i)$ starting at $c(t_i)$ converge to a
ray $c_\infty(t)$ in $X_\infty$ that starts at $p$, and the
corresponding Busemann functions also converge. Since the Busemann
functions on $X$ are $C^2$~\cite{HIH}, they converge to a $C^1$
Busemann function on $X_\infty$. Thus the horosphere passing
through $p$ is a $C^1$-submanifold of $X_\infty$, and is the limit
of horospheres passing through $c(t_i)$.
The Busemann flow is $C^1$ on $X$, and $C^0$ on $X_\infty$.
Since the horospheres in
$X$ and $X_\infty$ have the same dimension, the sequence of
horospheres passing through $c(t_i)$ does not collapse, and more
generally, each horosphere centered at $c_\infty(\infty)$ is the
limit of a noncollapsing sequence of horospheres in $X$.

It is easy to see that the group $G$ fixes $c_\infty(\infty)$, i.e.
any $\g\in G$ takes $c_\infty$ to a ray asymptotic to $c_\infty$.
Furthermore, $G$ leaves the horospheres corresponding to
$c_\infty(\infty)$ invariant.

Thus, one can define the rotation homomorphism $\phi_\infty\co
G\to O(n)$ corresponding to the point $c_\infty(\infty)$. The
point only determines $\phi_\infty$ up to conjugacy so we also
need to fix an isometry $L\co\Rn\to T_{p} X_\infty$. Similarly, a
choice of an isometry $L_i\co \Rn\to T_{c(t_i)} X$ specifies the
rotation homomorphism $\phi_i\co \G\to O(n)$ corresponding to the
point $c(\infty)$. We can assume that $\phi_i=\phi_0$ for each
$i$, by choosing $L_i$ equal to $L_0$ followed by the parallel
transport $P^\infty_{c(t_0),c(t_i)}=P_{c(t_0),c(t_i)}$. 
Henceforth we denote $\phi_0$ by $\phi$. Also it is
convenient to choose $L$ as follows.

\begin{lem} \label{lem: conv rotations}
After passing to a subsequence of $(X, c(t_i),\G)$,
there exist $L$ such that if $\g_i\to\g$, then
$\phi(\g_i)\to\phi_\infty(\g)$.
\end{lem}
\begin{proof}
Since $(X, c(t_i),\G)\to (X_\infty,p, G)$ in pointed equivariant
$C^{1,\a}$ topology, we can find the corresponding $C^{1,\a}$
approximations $f_i\co B(c(t_i),1)\to B(p,1)$. We may assume that
$df_i\co T_{c(t_i)} X\to T_p X_\infty$ is an isometry for all $i$.
By compactness of $O(n)$, $df_i\circ L_i$ subconverge to an
isometry $L\co\Rn\to T_p X$, so by modifying $f_i$ slightly we can
assume that $df_i\circ L_i=L$. This $L$ is then used to define
$\phi_\infty$, and it remains to show that if $\g_i\to\g$, then
$\phi(\g_i)\to\phi_\infty(\g)$. For the rest of the proof we
suppress $L_i, L$.

Since $d\gamma_i\to
d\gamma$ it is enough to show that parallel transports through
infinity from $c(t_i)$ to $\g_i(c(t_i))$ converge to the
parallel transport through infinity from $p$ to $\gamma(p)$.

For $x\in X$ and $c(t)$ that lie on the same horosphere, and for
any $s>t$, denote by $P_{x,s,c(t)}\co T_{c(t)}M\to T_x M$ the
parallel transport along the piecewise geodesic path $x ,x+s,
c(t)+s, c(t)$. Here $c(t)+s=c(t+s)$, $x+s$ also lie on the same
horosphere. Now $|P^\infty_{x c(t)}-P_{x,s,c(t)}|$ can be
estimated as
\[|P^\infty_{x+s,c(t)+s}-P_{x+s,c(t)+s}|\le q(a)d(x+s,c(t)+s)\le
q(a) e^{-s} c_2(a,d(x,c(t))).\] By Remark~\ref{rmk: valid for c1},
the same estimate holds for $X_\infty$ and $c_\infty$.

Fix $\e>0$ and pick $R>0$ such that $d(c(t_i),\g_i(c(t_i)))\le R$
for all $i$. Take a large enough $s$ so that $q(a)c_2(a, R) e^{-
s}<\e$. Since $B(c(t_i),R+s)$ converges to $B(p,R+s)$ in
$C^{1,\a}$-topology, and $\g_i\to \g$, we conclude that
$P_{\g_i(c(t_i)),s, c(t_i)}\to P_{\g(p),s,p}$ in $C^0$-topology,
or more formally,
\[|df_i\circ P_{\g_i(c(t_i)),s, c(t_i)}\circ
d\g_i-P_{\g(p),s,p}\circ d\g|<\e.\]
for large $i$, where $f_i$ is the $C^{1,\a}$-approximation.
By the estimate in the previous paragraph,
$|P_{\g_i(c(t_i)),s,c(t_i)}-P^\infty_{\g_i(c(t_i))c(t_i)}|<\e$
and $|P_{\g(p),s, p}-P^\infty_{\g(p),p}|<\e$,
so the triangle inequality implies that
$|df_i\circ P^\infty_{\g_i(c(t_i))c(t_i)}\circ d\g_i -
P^\infty_{\g(p),p}\circ d\g|<3\e$
for large $i$. Hence
$|\phi(\g_i)-\phi_\infty(\g)|< 3\e$ for all large $i$, and
as $\e>0$ is arbitrary it follows that
$\phi(\g_i)\conv{i\to\infty}\phi_\infty(\g)$.
\end{proof}

\begin{prop}\label{prop: G and K}
Let $K=\ker \phi_\infty$ and let $G_p$ be the isotropy subgroup of $p$ in $G$. Then\newline
$\mathrm(1)$ $\overline{\phi(\G)}=\phi_\infty(G_p)=\phi_\infty(G)$,\newline
$\mathrm(2)$  $K$ acts freely on $X$, in particular, 
$K\cap G_p=\{\mathrm{id}\}$.\newline
$\mathrm(3)$
The short exact sequence 
$1\to K\to G\overset{\phi_\infty}{\to}\phi_\infty(G)\to 1$ splits
with the splitting given by
$\phi_\infty(G)\simeq G_p\hookrightarrow G$.
In particular, $G$ is a semidirect product of $K$ and $G_p$.
\end{prop}
\begin{proof}
(1) For each $\g\in\G$ we have $d(\g(c(t_i)), c(t_i))\to 0$ as $i\to\infty$,
so the constant sequence $\g$ converges to some $g\in G_p$.
By Lemma~\ref{lem: conv rotations}, $\phi(\g)\to \phi_\infty(g)$,
which means $\phi(\g)=\phi_\infty(g)$. Thus,
$\phi(\G)\subset\phi_\infty(G_p)$.
Now $G_p$ is compact, so $\phi_\infty(G_p)$ is closed, and
therefore, $\overline{\phi(\G)}\subset\phi_\infty(G_p)$.
Since $\phi_\infty(G_p)\subset \phi_\infty(G)$, it remains to show
that $\phi_\infty(G)\subset\overline{\phi(\G)}$.
Given $\g\in G$, we find $\g_i\in\G$ with $\g_i\to\g$.
By Lemma~\ref{lem: conv rotations}, $\phi_\infty(\g)$ is the limit
of $\phi(\g_i)\in \phi(\G)$, so $\phi_\infty(g)\in\overline{\phi(\G)}$.

(2) If $k\in K$ fixes a point $x$,
then $1=\phi_\infty(k)=P^\infty_{k(x)x}\circ dk=P^\infty_{xx}\circ dk=dk$.
Since $k$ is an isometry, $k=\mathrm{id}$.
(3) is a formal consequence of (1) and (2).
\end{proof}

\begin{rmk}
Since $G$ is the semidirect product of $K$ and $G_p$,
any $\g\in G$ can be uniquely written as $kg$ with
$k\in K$, $g\in G_p$, We refer to $k$ and $g$ respectively as the
{\it translational part} and the {\it rotational part} of $\g$.
\end{rmk}

\begin{rmk} \label{rmk: G is almost nil}
If $\G$ is discrete, then by Margulis lemma any finitely generated
subgroup of $\G$ has a nilpotent subgroup whose index $i$ and the
degree of nilpotency $d$ are bounded above by a constant depending
only on $n$. (Of course, by~\cite{Bow} $\G$ itself is finitely
generated, as is any subgroup of $\G$, but we do not need this
harder fact here). The same then holds for $G$. Indeed, take
finitely many elements $g_l$ of $G$ and approximate them by
$\g_{j,l}\in\G$, so they generate a finitely generated subgroup of
$\G$. Then  $\g_{j,l}^i$ approximate $g_l^i$, and by above,
$g_l^i$ lie in a nilpotent subgroup of $\G$. Hence a $d$-fold
iterated commutator in $g_{j,l}^i$'s is trivial for all $j$, and
then so is the corresponding commutator in $g_l^i$'s. 
Hence $G$ is nilpotent by~\cite[Lemma VIII.8.17]{Rag}
\end{rmk}

\section{Controlling injectivity radius}

We continue working with notations of Section~\ref{sec: pass to
lim}, except now we also assume that $\G$ is discrete. The family
$(X, c(t),\G)$ may have many converging subsequences with limits
of the form $(X_\infty, p, G)$. We denote by $K(p)$ the $K$-orbit
of $p$, where $K$ is the kernel of the rotation homomorphism $G\to
O(n)$. The goal of this section is to find a common lower bound,
on the normal injectivity radii of $K(p)$'s.

\begin{prop}\label{2form}
There exists a constant $f(a)$ such that for each $x\in K(p)$, the
norm of the second fundamental form $II_x$ of $K(p)$ at $x$ is
bounded above by $f(a)$.
\end{prop}
\begin{proof}
Since $K$ acts by isometries $|II_x|=|II_p|$ for any $x\in K(p)$ so
we can assume $x=p$. Let $X,Y\in T_p K(p)$ be unit tangent
vectors. Extend $Y$ to a left invariant vector field on $K(p)$,
and let $\a(t)=exp(tX)(p)$ be the orbit of $p$ under the
one-parameter subgroup generated by $X$. Since $II_p(X,Y)$ is the
normal component of $\nabla_XY(p)$, it suffices to show that
$|\nabla_XY|\le f(a)$. Let $P^\a_{p,\a(t)}$ be the parallel
transport from $p$ to $\a(t)$ along $\a$. By Section~\ref{sec: rot
hom} $|P_{p,\a(t)}-P^\infty_{p,\a(t)}|\le q(a)d(p,\a(t))\le
2q(a)t$ for all small $t$.

A similar argument shows that $|P_{p,\a(t)}-P^\a_{p,\a(t)}|\le
2q(a) t$ for all small $t$. Indeed, look at the ``ruled'' surface
obtained by joining $p$ to the points of $\a$ near $p$. If we
approximate $\a$ by a piecewise geodesic curve
$p\a(t_1)\cdots\a(t_k)$, where $\a(t_k)=q$ is some fixed point
near $p$, then the area of the surface can be computed as the
limit as $k\to\infty$ of the sum of the areas of geodesic
triangles $p\a(t_i)\a(t_{i+1})$. The area of each triangle is
bounded above by $d(\a(t_i)\a(t_{i+1}))$, so the area of the ruled
surface is bounded above the length of $\a$ from $p$ to $q$, which
is at most $2t$, for small $t$.

Therefore, $|P^\infty_{p,\a(t)}-P^\a_{p,\a(t)}|\le 4q(a)t=f(a)t$
by the triangle inequality, so
$|P^\infty_{\a(t),p}Y-P^\a_{\a(t),p}Y|\le f(a)t$. On the other
hand, $P^\infty_{\a(t),p}Y=Y(p)$ because $Y$ is left-invariant,
and since elements of $K$ have trivial rotational parts. Thus
$|P^\a_{\a(t),p}Y-Y(p)|\le f(a) t$, which by definition of
covariant derivative implies that $|\nabla_XY(p)|\le f(a)$.
\end{proof}

\begin{cor} \label{cor: 2fund}
$\mathrm{(i)}$ There exists $r(a)>0$ such that if $C$ is the
connected component of $K(p)\bigcap B_{r(a)}(p)$ that contains
$p$, and if $x\in B_{r(a)}(p)$ is the endpoint of the geodesic
segment $[x,p]$ that is perpendicular to $C$ at $p$, then
$d(x,c)>d(x,p)$ for any $c\in C\setminus\{p\}$. \newline
$\mathrm{(ii)}$ If there exists $s<r(a)$ such that $K(p)\bigcap
B_{s}(p)$ is connected, then the normal injectivity radius of
$K(p)$ is $\ge s/3$.
\end{cor}
\begin{proof}
(i) The metric on $X_\infty$ can be approximated in $C^1$-topology
by smooth metrics with almost the same two-sided negative
curvature bounds and infinite injectivity radius~\cite{Nik}. Also
$C^1$-closeness of metrics implies $C^0$-closeness of Levi-Civita
connections, and hence almost the same bounds on the second
fundamental forms of $C$. Now for the smooth metrics as above the
assertion of (i) is well-known, and after choosing a slightly
smaller $r(a)$, it passes to the limits so we also get it for
$X_\infty$.

(ii) Consider two arbitrary geodesic segments of equal length $\le
s/3$ that start at $K(p)$, are normal to $K(p)$ and have the same
endpoint. Since $K$ acts isometrically on $X_\infty$ and
transitively on $K(p)$, we can assume that one of the segments
starts at $p$. By the triangle inequality, the other segment
starts at a point of $C$, so by part (i) the segments have to
coincide.
\end{proof}

\begin{rmk}
The proof that $K(p)\bigcap B_{s}(p)$ is connected for some
$s<r(a)$ independent of the converging subsequence
$(X,c(t_i),\Gamma)$ occupies the rest of this section, and this is
the only place in the paper that uses Bowditch's
theorem~\cite{Bow} that $\G$ is finitely generated. Other key
ingredients are the existence of approximate square roots in
finitely generated nilpotent groups (see Appendix~\ref{app: nilp}),
and the following comparison lemma that relates the
displacement of an element of $\G$ to the displacement of its
square root.
\end{rmk}

\begin{figure}[!h]
\caption{}
\includegraphics[scale=.75,trim=30mm 1mm 1mm 1mm]{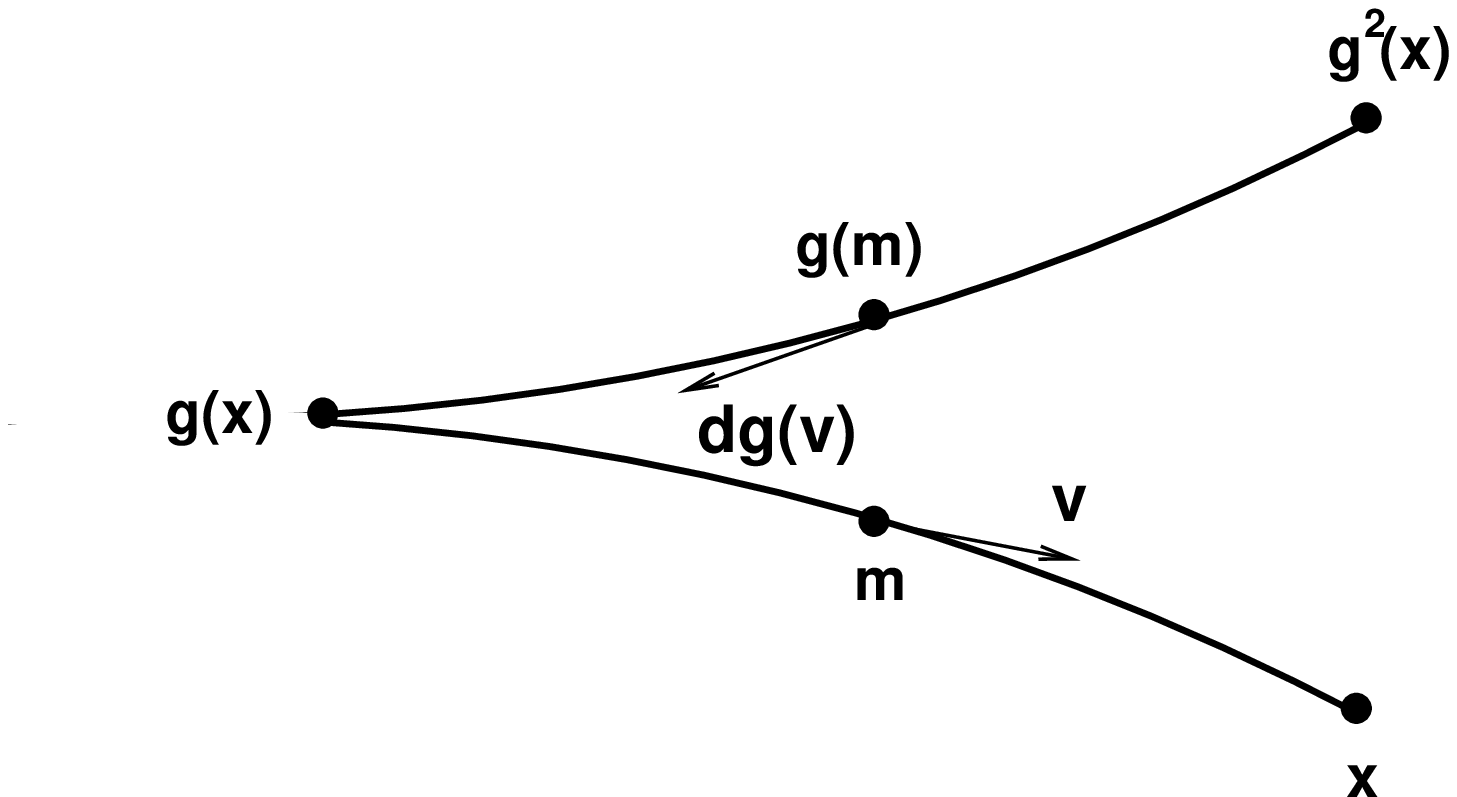}
\end{figure}

\begin{lem}\label{est}
Let $U$ be the neighborhood of $1\in O(n)$ that consists of all
$A\in O(n)$ satisfying $|Av-v|<1$ for any unit vector in
$v\in\Rn$. There exists a function $f\co (0,\infty)\to (0,\infty)$
such that $f(r)\to 0$ as $r\to 0$ and $d(g(x),x)\le
f(d(g^2(x),x))$, for any $x\in X$ and any $g\in \Gamma$ with
$\phi(g)\in U$.
\end{lem}
\begin{proof} Define $f(r)$ to be the supremum of numbers
$d(g(x),x)$ over all $x\in X$ and $g\in \Gamma$ with $\phi(g)\in
U$ satisfying $r=d(x,g^2(x))$.

To see that $f(r)<\infty$ take an arbitrary $x\in X$, $g\in\G$
with $r=d(x,g^2(x))$ and let $R=d(x,g(x))$. Look at the geodesic
triangle in $X$ with vertices $x$, $g(x)$, $g^2(x)$
(see Figure 1).
Arguing by
contradiction, assume that by choosing $g, x$ one can make $R$
arbitrary large while keeping $r$ fixed. The geodesic triangle
then becomes very long and thin. Let $m$ be the midpoint of the
geodesic segment $[x,g(x)]$, so that $g(m)$ is the midpoint of the
geodesic segment $[g(x),g^2(x)]$. By exponential convergence of
geodesics and comparison with the hyperbolic plane of $\sec=-1$,
we get $d(m,g(m))\le C(r) e^{-R/2}$, which is small since $R$ is
large. So $P^\infty_{g(m) m}$ is close to $P_{g(m) m}$, which in
turn is close to $P_{g(x) m}\circ P_{g(m) g(x)}$, since the
geodesic triangle with vertices $m$, $g(m)$, $g(x)$ has small
area. Thus, $\phi(g)$ is close to $P_{g(x) m}\circ P_{g(m)
g(x)}\circ dg$. Let $v$ be the unit vector tangent to $[x,g(x)]$
at $m$ and pointing towards $x$. Then $dg(v)$ is tangent to
$[g(x),g^2(x)]$ at $g(m)$ and is pointing towards $g(x)$. Since
the geodesic triangle with vertices $x$, $g(x)$, $g^2(x)$ has
small angle at $g(x)$, the map $P_{g(x) m}\circ P_{g(m) g(x)}$
takes $dg(v)$ to a vector that is close to $-v$. This gives a
contradiction since $|\phi(g)(v)-v|\le 1$.

A similar argument yields that $f(r)\to 0$ as $r\to 0$. Namely, if
one can make $d(g^2(x),x)$ arbitrary small while keeping
$d(g(x),x)$ bounded below, then the geodesic triangle with
vertices $x$, $g(x)$, $g^2(x)$ becomes thin, and we get a
contradiction exactly as above.
\end{proof}

\begin{prop}\label{dist}
Let $r(a)$ be the constant of Corollary~\ref{cor: 2fund}. Then
there exists a positive $s\le r(a)$, depending only on $X$, $c$,
$\G$, such that for any converging sequence $(X, c(t_i), \G)\to
(X_\infty,p,G)$, the normal injectivity radius of $K(p)$ is $\ge
s$.
\end{prop}
\begin{proof}
By Corollary~\ref{cor: 2fund}, it suffices to find a universal $s$
such that $K^s:=K(p)\bigcap B_{s}(p)$ is connected. Let $K^s_0$ be
the component of $K^s$ containing $p$.

By~\cite{Bow} $\G$ is finitely generated, hence by Margulis
lemma~\cite{BGS}, $\G$ contains a normal nilpotent subgroup
$\tilde{\G}$ of index $i\le i(a,n)$. Therefore
$H=\overline{\phi(\G)}$ is virtually nilpotent. Hence, its
identity component $H_0$ is abelian, since compact connected
nilpotent Lie groups are abelian. Then $\phi^{-1}(H_0)$ is a
subgroup of finite index in $\Gamma$. Let $\G'=\tilde{\G}\cap
\phi^{-1}(H_0) $. Clearly, $[\Gamma:\G^\prime]= k=k(\phi )$ 
is also finite.

We first give a proof in case $\G=\G^\prime$.
Arguing by contradiction, suppose that  for any $s>0$ there exists
a sequence $(X, c(t_i), \Gamma)\to (X_\infty,p,G)$ and a point
$\g(p)\in K^s\setminus K^s_0$ with $d(p,\g(p))<s$. By possibly
making $d(p,\g(p))$ smaller, we can choose $\g(p)$ so that
$d(p,\g(p))$ is the distance from $p$ to $K^s\setminus K^s_0$. By
the first variation formula the geodesic segment $[p,\g(p)]$ is
perpendicular to $K^s_0$. The next goal is to construct the square
root of $\g$ with no rotational part and displacement bounded by
$f(d(p,\g(p)))$, where $f$ is the function of Lemma~\ref{est}.

Take $\g_i\in \G$ that converge to $\g$. Since $\G$ is finitely
generated, we can apply Lemma~\ref{root} to find an (independent of
$i$) finite set $F\subset \G$ such that each $\g_i$ can be written
as $\g_i=g_i^2f_i$ with $f_i\in F$. We can further write each
$g_i$ as the product $g_i=x_i r_i$, where $x_i$ has a small
rotational part and $r_i$ is close to a rotation, namely we let
$r_i$ be an element of $\G$ that is close to $\phi(g_i)$ and let
$x_i=g_i r_i^{-1}$. Thus $\g_i=(x_i r_i)^{2}f_i$ and
\[
\g_i=x_i r_i x_i r_i f_i= x_i^{2} [x_i^{-1} r_i] r_i^2 f_i
\]
Applying Lemma~\ref{root2} to $x^2_i[x_i^{-1} r_i]$ we see that
$\g_i$ can be written as $(x_i h_i)^{2}f_{i}^\prime r_i^{2}f_i$
with $h_i\in [\G,\G]$, $f_i^\prime\in F^\prime$. Since $\phi(\G)$
ia abelian, we have $\phi(h_i)=1$ and hence $\phi(x_i
h_i)=\phi(x_i)$ is small. Since $F$, $F^\prime$ are independent of
$i$, each element of $F$, $F^\prime$ is close to a rotation for
large $i$. So $f_{i}^\prime r_i^{2}f_i$ is close to a rotation and
since both $\phi(\g_i)$ and $\phi(x_ih_i)$ are small,
$f_{i}^\prime r_i^{2}f_i$ subconverges to the identity. Thus, $(x_i
h_i)^{2}$ subconverges to $\g$, and so we might as well assume
that $\gamma_i=(x_ih_i)^2$ in the beginning. By Lemma~\ref{est},
\[
d(c(t_i),x_i h_i(c(t_i)))\le f(d(c(t_i),(x_i h_i)^2(c(t_i)))
\]
where the right hand converges to $f(d(p,\g(p)))$. Hence $x_i h_i$
subconverges to $w\in K$ such that $w^2=\g$ and $d(p,w(p))\le
f(d(p,\g(p)))$.

Since $w$ has no rotational part, $P^\infty_{pw(p)}=dw_p$. By
assumption $s$ can be taken arbitrary small, so we can assume that
$f(d(p,\g(p)))$ is small, in particular, $w(p)\in K^s$. So
$P_{pw(p)}$ is close to $dw_p$. Hence if $v$ is the unit vector
tangent to $[p,w(p)]$ at $p$ and pointing towards $w(p)$, then
$P_{pw(p)}(v)$ is close to $dw_p(v)$. Therefore, $w(p)$ is close
to the midpoint of $[p,w^2(p)]=[p,\g(p)]$. Since $[p,\g(p)]$ is
perpendicular to $K^s_0$ and $d(p,\g(p))<r(a)$, it is clear that
$w(p)\notin K^s_0$. This contradicts the minimality of
$d(p,\g(p))$ and completes the proof in case $\G=\G^\prime$.

We now turn to the general case.  Let $G^\prime$ be the subset of
$G$ that consists of limits of elements of $\G^\prime$ under the
convergence $(X, c(t_i), \G)\to (X_\infty,p,G)$. It is
straightforward to check that $G^\prime$ is a closed subgroup of
$G$ of index $\le k$. Thus the limit of any converging subsequence
of $(X, c(t_i), \G^\prime)$ has to equal to
$(X_\infty,p,G^\prime)$, therefore in fact, $(X, c(t_i),
\G^\prime)$ converges to $(X_\infty,p,G^\prime)$. Since the
rotation homomorphism of $G$ restricts to the rotation
homomorphism of $G^\prime$, the translational part $K^\prime$ of
$G^\prime$ is $G^\prime\cap K$. In particular, $|K:K^\prime|\le k$
hence the identity components of $K$ and $K^\prime$ coincide.
Using the first part of the proof, we fix $s$ such that
$K^\prime(p)\bigcap B_s(p)$ is connected. Thus $K^\prime(p)\bigcap
B_s(p)=K^s_0$.

Now let $\g(p)\in K^s\setminus K^s_0$ such that $d(p,\g(p))$ is
the distance from $p$ to $K^s\setminus K^s_0$. Then the geodesic
segment $[p,\g(p)]$ is perpendicular to $K^s_0$ by the first
variation formula. Arguing by contradiction suppose that
$d(p,\g(p))$ can be arbitrary small. Then by triangle inequality
$\g^j(p)$ is close to $p$ for $j=1,\dots, k$. Also $\g^k\in
K^\prime$, so in fact $\g^k(p)\in K^s_0$ because
$K^\prime(p)\bigcap B_s(p)=K^s_0$.

On the other hand, since $\g^j$'s have no rotational part, the
argument used above to prove that $w(p)$ is close to the midpoint
of $[p,w^2(p)]$ shows that the points $\g^j(p)$ almost lie on a
geodesic segment $[p, \g^k(p)]$. Then the segments $[p, \g^k(p)]$
and $[p,\g(p)]$ have almost the same direction, so $[p, \g^k(p)]$
is almost perpendicular to $K^s_0$. Hence by Corollary~\ref{cor:
2fund}, if $d(p,\g(p))$ is small enough, then $\g^k(p)\notin
K^s_0$ which is a contradiction.
\end{proof}

\begin{rmk}
Although it is not needed for the proof of Theorem~\ref{main thm},
note that all possible limits of $(X/\G,\s(t_i))$ have the
same dimension independent of the sequence $t_i\to\infty$.
Consider all possible limits with fixed $\dim(K)$, and look at the space
$X_\infty/ K$. It has a lower bound on the injectivity radius for
points near the projection $\bar p$ of $p$,
and hence $\vol(B(\bar p,1))>c>0$
in all such spaces, where $c$ is independent
of the converging sequence.
By Proposition~\ref{prop: G and K},
the isotropy group $G_p$ is the same for all possible limits
and moreover, the $G_p$-actions on $T_pX_\infty$ are all equivalent.
Also by Lemma~\ref{lem: max comp} the identity component
$G_p^\mathrm{id}$ of $G_{p}$
commutes with the identity component of $K$, hence
$G_p^\mathrm{id}$ fixes pointwise
the component of $K(p)$ containing $p$. Thus
the $G_p^\mathrm{id}$-actions
on $T_{\bar p}X_\infty/ K$ are all equivalent.
This implies that a unit ball in $X_\infty/G$ has
volume $>c^\prime>0$ with $c^\prime$ only depending
on $\dim (K)$.
Hence all limits of the same dimension form a closed subset among
all limits, therefore, the space of all limits is the union of
these closed sets. On the other hand, the space of all limits is connected
by Lemma~\ref{lem: connected} below,
thus all the limits have the same dimension.

\begin{lem} \label{lem: connected}
If $\g\co [0,\infty)\to Z$ is a continuous precompact
curve in a  metric space $Z$, then the space $Lim(\g)$ of all
possible subsequential limits $\lim_{t_i\to\infty}\g(t_i)$ is
connected.
\end{lem}
\begin{proof}
If $Lim(\g)$ is not connected, then we can write it as a disjoint
union of closed (and hence compact) sets $Lim (\g)=A\sqcup B$.
Then $U_\e(A)\cap U_\e(B)=\emptyset$ for some $\e>0$,
where $U_\e (S)$ denotes the $\e$-neighborhood of $S$. Let
$\g(t_i)\to a\in A$ and $\g(t_i')\to b\in B$. Arguing by
contradiction, we see that the curve $\g|_{[t_i,t'_i]}$ lies in
$U_\e(Lim(\g))=U_\e(A)\cup U_\e(B)$ for all large $i$. Clearly,
$\g(t_i)\in U_\e(A)$ and $\g(t_i')\in U_\e(B)$ for all large $i$,
which contradicts $U_\e(A)\cap U_\e(B)=\emptyset$.
\end{proof}
\end{rmk}

\section{Product structure at infinity}\label{prod}

In the next two sections we apply the critical point theory for
distance functions to show the following.

\begin{thm}\label{thm: prod at infty}
For each large $t$, the horosphere quotient $H_{t}/\G$ is
diffeomorphic to the normal bundle of an orbit of an N-structure
on $H_{t}/\G$.
\end{thm}
\begin{proof}
Let $\s(t)$ be the projection of $c(t)$ to $X/\G$. Given a
converging sequence $(X, c(t_i),\G)\to (X_\infty,p,G)$, the
sequence of pointed Riemannian manifolds $(X/\G,\s(t_i))$
converges in the pointed Gromov-Hausdorff topology to a pointed
Alexandrov space $(Y,q):= (X_\infty/G,q)$ with curvature bounded below
by $-a^2$.

The identity component $K_\mathrm{id}$ of $K$ is normal in $K$,
hence its $p$-orbit $K_\mathrm{id}(p)$ is invariant under the
action of $G_p$. Let $r\ll s$ be a positive constant to be
determined later, where $s$ comes from Proposition~\ref{dist}. The
$3r$-tubular neighborhood of $K_\mathrm{id}(p)$ is also
$G_p$-invariant, so the ball $B_{3r}(q)$ is isometric to the
$G$-quotient of this tubular neighborhood. By
Proposition~\ref{dist}, any $x\in B_{3r}(q)$ can be joined to $q$
by a unique shortest geodesic segment $[q,x]\subset B_{3r}(q)$.

Recall that in general a distance function $d(\cdot,q)$ on an
Alexandrov space is called {\it regular at the point $x$} if there
exist a segment emanating from $x$ that forms the angle $>\pi/2$
with any shortest segment joining $x$ to $q$.

In our case the function $d(\cdot ,q)$ is
regular at any $x\in B_{3r}(q)\setminus\{q\}$.
%

Let $w(a)>1$ be a constant depending only on $a$ that will be
specified later. By angle comparison the function $d(\cdot,\s(t))$
on
\[A_r(\s(t))=\{x\in B_s(\s(t))\co d(x,\s(t))\in [r/w(a),w(a)r]\}\]
is regular provided the Gromov-Hausdorff distance between $B_s(q)$
and $B_s(\s(t))$ is $\ll r/w(a)$. Because the family
$\{B_s(\s(t))\}$ is precompact in the Gromov-Hausdorff topology, 
Proposition~\ref{dist} implies that 
the function $d(\cdot ,q)$ is regular on
$A_r(\s(t))$ for all $t\ge t_0$ with sufficiently large $t_0$.

We denote by $H_t$ the horosphere centered at $c(\infty)$ that
contains $c(t)$. Since the second fundamental form of $H_t$, is
bounded in terms of $a$, any short segment joining  nearby points of
$H_t/\G$ is almost tangent to $H_t/\G$. Hence by taking $r$
sufficiently small, we can assume that for all $t\ge t_0$ and each
$x\in A_r(\s(t))$ there exists a unit vector $\l\in T_x (H_t/\G)$
that forms the angle $\a_{\l, [x,\sigma(t)]}\in
[\frac{2\pi}{3},\pi]$ with any shortest segment $[x,\sigma(t)]$.
By the first variation formula, the derivative of $d(\cdot ,q)$ in
the direction of $\l$ equals to the minimum of $-\cos\a_{\l,
[x,\sigma(t)]}$ over all shortest segments $[x,\sigma(t)]$ and by
above it  lies in $ [\frac{1}{2}, 1]$.

The distance function on $X/\G$ need not be smooth, and for what
follows it is convenient to replace $d(\cdot, \s(t))$ by its
average over a small ball $B_\d(\s(t))$ as follows. Given
$\delta\ll r$, define $f\co X/\G\to\mathbb R$ by
\[f(x)=\frac{1}{\vol B_\d(\s(t))}\int_{B_\d(\s(t))} d(x,y)dy\]
where $x\in H_t/\G$(i.e. $t=b(x)$). Now $f$ is $C^1$ $1$-Lipschitz
function with $|f(x)-d(x,\s(t))|\le\d$ for any $x\in A_r(\s(t))$.
Observe that for any $\eta\in T_x(X/\G)$
\begin{equation}\label{der1}
df_x(\eta)=\frac{1}{\vol B_\d(\s(t))}\int_{B_\d(\s(t))}
(-\cos\a_{\eta,[x,y]})dy
\end{equation}

Also note that since $\d \ll r$ and $\sec(X/\G)\ge -a^2$, for all large
$t$ if $x\in A_r(\s(t))$ and $y\in B_\d(\s(t))$ there is a point
$z$ such that $d(z,x)\approx d(x,y)$ and $d(z,y)\approx 2d(x,y)$.
Therefore, the
 angle corresponding to $x$ in the comparison triangle in the space of
$\sec\equiv -a^2$ is almost $\pi$. By Toponogov comparison, the
angle at $x$ in any geodesic triangle $\Delta xyz$ is almost
$\pi$.  By~(\ref{der1}), this implies that if $\eta$ is a direction
of any shortest segment connecting $x$ to $z$, then $df_x(\eta)\in
[\frac{1}{2}, 1]$ provided $\d$ is small enough.

By gluing $\l$'s via a partition of unity, we obtain a $C^1$ unit
vector field $\L$ that is tangent to $H_t/\G$ and defined for all
$t\ge t_0$ and $x\in A_r(\s(t))$, and such that $df_x(\L)\in
[\frac{1}{4}, 1]$ if $\d$ is sufficiently small. Then $S_r=\{x\in
X/\G: f(x)=r\}$ is a properly embedded $C^1$-hypersurface in
$X/\G$ that is transverse to $\L$. Also the compact submanifold
$S_t(t):=S_r\bigcap H_t/\G$ is $\d$-close to the metric $r$-sphere
in $H_t/\G$ centered at $\s(t)$. Furthermore,
$A(r,t):=A_r(\s(t))\cap H_t/\G$ is $C^1$-diffeomorphic to the
product $S_r(t)\times [r/w(a) ,w(a)r]$.

Here we are only interested in the part of $X/\G$ with $t\ge t_0$.
There the Busemann function $X/\G\to\mathbb R$ restricts to a
$C^1$-submersion $S_r\to\mathbb [t_0,\infty)$, because otherwise
at some point the tangent spaces of $S_r$ and $H_t/\G$ coincide by
dimension reasons, so that $S_r$ cannot be transverse to $\L\in
T(H_t/\G )$. By construction the submersion is proper, hence it is a
$C^1$-fiber bundle, which is $C^1$-trivial by the covering
homotopy theorem. The trivialization defines a $C^1$-isotopy $F\co
S_r(t_0)\times [t_0,\infty)\to X/\G$ such that
$S_r(t_0)\times\{t\}$ is mapped onto $S_r(t)$.

We push this isotopy along the Busemann flow back into
$H_{t_0}/\G$ by setting $G(t,x)=b_{t_0-t}(F(x,t))$ for any $t\ge
t_0$, and $x\in S_r(t_0)$, to get the $C^1$-isotopy $G\co
S_r(t_0)\times [t_0,\infty)\to H_{t_0}/\G$.

The Busemann flow induces a $C^1$-diffeomorphism $H_t/\G\to
H_{t_0}/\G$, so around each submanifold $b_{t_0-t}(S_r(t))$ there
is a ``tubular neighborhood'' $b_{t_0-t}(A(r,t))$.

By the exponential convergence of geodesics, one can choose $w(a)$
in the definition of $A_r(\s(t))$ so that for any $t\ge t_0$ there
exists $t^\prime\ge t+1$ such that
$b_{t_0-t^\prime}(S_r(t^\prime))$ is contained in
$b_{t_0-t}(A(r,t))$ and is disjoint from $b_{t_0-t}(S_r(t))$. By
the following elementary lemma, the region between
$b_{t_0-t}(S_r(t))$ and $b_{t_0-t^\prime}(S_r(t^\prime))$ is
$C^1$-diffeomorphic to $S_r(t)\times [0,1]$.

\begin{lem}\label{l:isotopy}
Let $M$ a closed smooth manifold $M$ and $F_t\co M\to
M\times\mathbb R$ be a $C^1$-isotopy with $F_0(M)=M\times\{0\}$.
If $F_0(M)$ and $F_s(M)$ are disjoint for some $s$, then the
region between $F_0(M)$ and $F_s(M)$ is diffeomorphic to $M\times
[0,1]$.
\end{lem}
\begin{proof}
By the isotopy extension theorem~\cite[p 293]{Cerf}
we can extend the isotopy $F_t$ to an ambient $C^1$-isotopy which
is identity outside a compact subset of $M\times\mathbb R$. Assume
without loss of generality that $s<0$, and then take $n>0$ so
large that the isotopy is identity on $M\times\{n\}$. Then by
restricting the ambient isotopy to the region between $M\times
\{n\}$, $M\times \{0\}$, we get a diffeomorphism of the region
between the region between $M\times \{n\}$, $M\times \{0\}$ onto
the region between $M\times \{n\}$, $F_s(M)$. The former region is
the product, so is the latter. But the latter region is
diffeomorphic to the region between $M\times \{0\}$, $F_s(M)$,
because the region between $M\times \{n\}$, $M\times \{0\}$ is
$M\times [0,1]$.
\end{proof}

By gluing a countable number of such diffeomorphisms together we
conclude that for all sufficiently large $t$
\begin{equation}\label{prodinf}
\begin{array}{c}
(H_{t}/\G)\setminus U(r,t) \text{ is $C^1$ diffeomorphic to }
[t,\infty )\times S_r(t) \\ \text{ where } U(r,t)=\{x\in
H_{t}/\G\co f(x)<r\}
\end{array}
\end{equation}

\section{Tubular neighborhood of an orbit}
\label{sec: tub nbhd}

It remains to understand the topology of $U(r,t)$, and we do so
for large enough $t$ and small enough $r$. The proof involves
the collapsing theory developed in~\cite{CFG} and the geometry of
Alexandrov spaces (for which we refer to~\cite{BGP} and 
Appendix~\ref{A:crit}).

Let us look at a converging sequence $(H_{t_i}/\G,\s(t_i))\to(H,q)$.
First, we replace the metric on $H_{t_i}/\G$ with an invariant
Riemannian metric which is $\e_i$ close to $H_{t_i}/\G$ in $C^1$
topology and  is $A(\e_i)$-regular~\cite{Shi,CFG, Nik}, where
$\e_i\to 0$ as $i\to\infty$. Also,  spaces $H_{t_i}/\G$ with the new
metrics have uniform curvature bound $|\sec|\le
C'=C'$~\cite{Shi} that depends only on the original curvature bound
of $H_{t_i}/\G$. The collapsing theory~\cite{CFG} yields, for each $i$, 
the following commutative diagram given by the invariant metric 
$h_i$ on $H_{t_i}/\G$.  
\begin{equation*}
\xymatrix{&FB_i\ar[d]\ar[r]^{\eta_i}&Y_{i}\ar[d]\\
&B_i\ar[r]^{\bar{\eta}_i} &X_{i}}
\end{equation*}
Here $B_i$ is the ball $B(\s(t_i),1)$ in the metric $h_i$, and $FB_i$
is the frame bundle of $B_i$.
The vertical arrows are quotient maps under isometric
$O(n)$-actions and $\eta_i$ is a Riemannian submersion given by
the $N$-structure on $FB_i$. Clearly the induced map
$\bar{\eta}_i$ is a submetry (see Appendix~\ref{A:crit} for
background on submetries). By  Lemma~\ref{sub}, the
Toponogov comparison  with $\curv\ge -C'$  holds for 
any triangle with vertices in 
$B(\bar{\eta}_i(\s(t_i)),\frac{1}{4})$ for all large $i$.

Since we will only be interested in the geometry of
$X_i$ inside $\frac{1}{8}$-neighborhood of $\bar{\eta}_i(\s(t_i))$, 
we will treat $X_i$'s as Alexandrov spaces.

Note that
$X_i\overset{G-H}{\longrightarrow} \bar{X}=B(q,1)$ 
and $\dim X_i=\dim \bar{X}$ for all large $i$.
We claim that there  exists an $R>0$  
and a sequence $q_i\in X_i$ converging to $q$ 
such that $d(\cdot, q_i)$ has no
critical points in $B(q_i,R)$ for all large $i$.
Consider two cases depending on whether $q$ lies 
on the boundary of the Alexandrov space $\bar{X}$:

{\bf Case 1:} Suppose $q\notin \partial \bar{X}$.

Since $\bar{X}$ has $\curv\ge -C'$ in comparison sense, by~\cite{Kap},
the exists a strictly concave function $u$ on a $B(q,R)$ for some
$R\ll 1$ such that it has a maximum at $q$ and the superlevel sets are
compact.  By possibly making $R$ smaller we can assume that 
$R<d(q,\partial \bar{X})$.
This function is constructed by taking averages and
minima of distance functions. Therefore, it naturally lifts to a
function $u_i$ on $X_i$ such that $u_i$ converges uniformly to $u$.
By~\cite[Lemma 4.2]{Kap}, the lifts $u_i$ are strictly concave on
$B(\bar{\eta}_i(\s(t_i)),R/2)$ for all large $i$. Let $q_i$ be the point of
maximim of $u_i$. By uniqueness of the maximim $q_i\conv
{i\to\infty}q$. By Lemma~\ref{crit}, $d(\cdot, q_i)$ has no
critical points in $B(q_i,R/3)$ for all large $i$.

{\bf Case 2:} Suppose now that  $q\in \partial \bar{X}$.  
Denote by $D\bar{X}$ and $DX_i$ the doubles of $\bar{X}$ and
$X_i$ along the boundary and let $\iota$ be 
the canonical involution.  By~\cite{Per}, the doubles are 
also Alexandrov spaces with $\curv \ge -C'$. 
It is clear that $DX_i\overset{G-H}{\longrightarrow} D\bar{X}$.
By construction, we can chose $u$ and $u_i$ to be $\iota$-invariant.  
As before, let $q_i$ be the point of
maximum of $u_i$.
Since $u_i(q_i)=u_i(\iota(q_i))$, 
by uniqueness of maximums of  strictly concave functions 
we see that $q_i$ must lie on $\partial X_i$.  
Again by Lemma~\ref{crit},  $d(\cdot, q_i)$ has no
critical points in $B(q_i,R/3)$  in $DX_i$ for 
all large $i$ and hence the same is true for $d(\cdot, q_i)$ in 
$B(q_i,R/3)$.
 
This immediately implies
that the distance function $d(\cdot, O_i)$ to the orbit
$O_i$ over $q_i$ has no critical points in the $R/3$ neighborhood
$U_{R/3}(O_i)$ of $O_i$ for all large $i$.
Indeed, let $x\in U_{R/3}(O_i)\backslash O_i$ and let $\g(t)$ be a geodesic
starting at $\bar{\eta}_i(x)$ such that
$\frac{d}{dt}d(\g(\cdot),q_i)'|_{t=0}>0$.
Since $\bar{\eta_i}\co (H(t_i)/\G,h_i)\to X_i$ is a submetry,
there exists $\tilde{\gamma}$,
a horizontal lift of $\g$ starting at $x$. Then
$d(\tilde{\g}(t),O_i)=d(\g(t),q_i)$ and hence
$d(\tilde{\g}(t),O_i)'|_{t=0}=d(\g(t),q_i)'|_{t=0}>0$.

Therefore $U_{r}(O_i)$ is
diffeomorphic to the total space of the normal bundle to $O_i$ in
$H_{t_i}/\G$ for any $r\le R/3$.
Since $O_i$ is Hausdorff close to $\s(t_i)$, the same is true for
$U(r,t_i)$.

Combining this with (\ref{prodinf}),
we conclude that $H_{t_i}/\G$ is diffeomorphic to
the total space of the normal bundle to $O_i$ in
$H_{t_i}/\G$ for all sufficiently large $i$.
Finally, since the above proof works for any sequence $t_i\to\infty$,
arguing by contradiction we conclude that $H_t/\G$ is
diffeomorphic to the normal bundle of an orbit of an $N$-structure
for all sufficiently large $t$.
This completes the proof of Theorem~\ref{thm: prod at infty}.
\end{proof}

\begin{rmk}
The reader may be wondering why we work with
the Alexandrov spaces $X_i$ instead of the Riemannian
manifolds $Y_i$.
This is because the curvature of $Y_i$ may tend to $\pm\infty$
as $i\to\infty$, which makes it hard to control the geometry
of $Y_i$'s. If instead of  $\e_i\to 0$, we take
$\e_i$  equal to a small positive constant $\e$, then $|\sec(Y_i)|\le C(\e)$,
but then it may happen that the injectivity radius of the
GH-limit of $Y_i$'s is $\ll\e$, so we cannot
translate the lower bound on the injectivity radius
from $Y_i$ to $H_{t_i}/\G$.
\end{rmk}

\begin{rmk} By Theorem~\ref{thm: prod at infty}
each orbit $O_{q_i}$ as above is homotopy equivalent to $X/\G$.
Thus all $O_{q_i}$'s are homotopy equivalent, hence they are all
affine diffeomorphic (see e.g.~\cite[Theorem 2]{Wil2}).
\end{rmk}

\section{Normal bundle is flat}

\begin{thm}\label{thm: norm bundle flat}
For each large $t$, the horosphere quotient $H_{t}/\G$ admits
an $N$-structure that has an orbit $O_t$ such that the normal bundle
to $O_t$ is a flat Euclidean vector bundle with total space
diffeomorphic to $H_{t}/\G$.
\end{thm}

Arguing by contradiction, it suffices to prove the theorem
for any sequence $t_i\to\infty$ such that
$H_{t_i}/\G$ converges in pointed Gromov-Hausdorff topology.
We fix such a sequence and assume for the rest of the proof
that $t$ belongs to the sequence.

We denote by $g_t$ the $C^1$-Riemannian metric on the horosphere
quotient $H_t/\Gamma$ induced by the ambient metric $(M,g)$.
Fix a small positive $\e>0$ to be determined later; this constant
will only depend on $(M,g)$. By Theorem~\ref{thm: prod at infty}
and~\cite{CFG, Nik}, for all large $t$ there exists an N-structure
on $H_t/\G$ with an orbit $O_t$ such that the normal bundle to
$O_t$ is diffeomorphic to $H_{t}/\G$. Also all orbits of the
N-structure have diameter $<\e$ with respect to an invariant
metric $h_t$ that is $\e$-close to $g_t$ in uniform
$C^1$-topology, i.e. $|g_t-h_t|<\e$,
$|\nabla_{g_t}-\nabla_{h_t}|<\e$. It remains to show that for all
large $t$, the normal bundle to $O_t$ in $H_t/\Gamma$ is flat
Euclidean. The proof break into two independent parts.

In Section~\ref{subsec: norm bundle in strat} we find a stratum
$F_t$ of the N-structure on $H_t/\Gamma$ such that $F_t$ is an
$h_t$-totally-geodesic closed submanifold that contains $O_t$ with
flat normal bundle. This uses only general properties of
N-structures.

In Section~\ref{subsec: norm bundle to strat} we show that the
restriction to $O_t$ of the normal bundle of $F_t$ in $H_t/\Gamma$
is flat, for all large $t$. This uses the flat connection of
Section~\ref{sec: rot hom} and the fact that $O_t\to H_t/\Gamma$
is a homotopy equivalence.
\subsection{Normal bundle in a stratum is flat}
\label{subsec: norm bundle in strat} Throughout
Section~\ref{subsec: norm bundle in strat} we suppress the index
$t$, and write $O$ in place of $O_t$, etc. Let $V$ be the tubular
neighborhood of $O$ that is sufficiently small so that all orbits
in $V$ have dimension $\ge\dim(O)$. Let $\tilde O$, $\tilde V$ be
their universal covers. According to~\cite[pp. 364--365]{CFG}, the
group $\Iso(\tilde V)$ contains a connected (but not necessarily
simply-connected) nilpotent subgroup $N$ that stabilizes $\tilde
O$, acts transitively on $\tilde O$, and also
$\Lambda=N\cap\pi_1(V)$ is a finite index subgroup of $\pi_1(V)$
and is a lattice in $N$. The following lemma is implicit
in~\cite{CFG}.

\begin{lem}
The subgroup $H$ of $\Iso(\tilde V)$ generated by $N$ and
$\pi_1(V)$ is closed, $N$ is the identity component in $H$ and the
index of $N$ in $H$ is finite.
\end{lem}
\begin{proof}
$\Lambda$ is a cocompact discrete subgroup of $N$ and also of its
closure $\bar N$ in $\Iso(\tilde V)$. Since $\dim(N)$, $\dim(\bar
N)$ are both equal to the cohomological dimension of $\Lambda$, we
get $N=\bar N$. Now let $\Lambda_0$ be a maximal finite index {\it
normal} subgroup of $\pi_1(\tilde V)$ that is contained in
$\Lambda$. If $\gamma\in\pi_1(V)$, then $N\cap \gamma
N\gamma^{-1}$ contains $\Lambda_0$ as a cocompact discrete
subgroup, so as before $\dim(N)$, $\dim(N\cap \gamma
N\gamma^{-1})$ are both equal to the cohomological dimension of
$\Lambda_0$, so $N=N\cap \gamma N\gamma^{-1}$, and $N$ is
normalized by $\pi_1(V)$. Thus $N$ is normal in $H$, and
$\Lambda=\Lambda_0$. Since $N$ is connected, it remains to show
that $|H:N|$ is finite. Since $N$ and $\pi_1(V)$ generate $H$, the
finite subgroup $\pi_1(V)/\Lambda$ of $H/N$ generates $H/N$, hence
$\pi_1(V)/\Lambda=H/N$.
\end{proof}

Let $I$ be the intersection of the isotropy subgroups of $H$ of
the points of $\tilde O$. Since $\tilde O$ is $H$-invariant, $I$
is normal in $H$. The fixed point set of $I$ is a totally geodesic
submanifold $\tilde F$ of $\tilde V$. The $H$-action on $\tilde F$
descends to an $H/I$-action on $\tilde F$. Since $\pi_1(V)$ is
torsion free and discrete, $\pi_1(V)\cap I$ is trivial, and we
identify $\pi_1(V)$ with its image in $H/I$. Denote the projection
of $\tilde F$ into $V$ by $F$.

\begin{lem}
The normal bundle to $O$ in $F$ is flat.
\end{lem}
\begin{proof}
The group $N$ acts transitively on $\tilde O$ so all isotropy
subgroups for the $N$-action on $\tilde O$ are conjugate. Since
they are also compact, they lie in the center of $N$~\cite{Iwa},
in particular all the isotropy subgroups are equal, and hence each
of them is equal to $I\cap N$. In particular, $N/(I\cap N)$ acts
freely and transitively on $\tilde O$. Since $\tilde O$ is
simply-connected, so is $N/(I\cap N)$. Thus, $I\cap N$ is the
maximal compact subgroup of $N$, hence by Lemma~\ref{lem: max
comp} $I\cap N$ is a torus, which we denote $T$. The torus is the
identity component of the compact group $I$, because $|I:I\cap
N|\le |H:N|<\infty$. Since $N/T$ acts freely and transitively on
$\tilde O$, we can choose a trivialization of $\nu$, the normal
bundle of $F$ in $\tilde V$ that is invariant under the left
translations by $N/T$. Namely, let $e\in\tilde O$ be the point
corresponding to $1\in N/T$ under the diffeomorphism
$N/T\cong\tilde O$. Fix an isomorphism $\phi\co\nu_e\to\{e\}\times
\mathbb R^k$, and then extend it to the $N/T$-left-invariant
isomorphism $\nu\cong \tilde O\times\mathbb R^k$. Now take
$\gamma\in\pi_1(V)$, and $x\in \tilde O$. Using the above
trivialization we define the {\it rotational part} of $\gamma$ as
the automorphism of $\{e\}\times\mathbb R^k$ given by
\[\phi\circ dL_{\gamma(x)^{-1}}\circ d\gamma\circ dL_x\circ\phi^{-1},\]
where $dL_x$ is the differential of the left translation by $x\in
N/T$. Since $O$ is an infranilmanifold, $\pi_1(V)$ acts on $\tilde
O$ by affine transformations, that is if $y\in \tilde O$, then
$\gamma(y)=n_\gamma\cdot A_\gamma(y)$, where $n_\gamma\in N/T$ and
$A_\gamma$ is a Lie group automorphism of $N/T$. Hence, for $z\in
\tilde O$, we get
\[(L_{\gamma(x)^{-1}}\circ \gamma\circ L_x) (z)=
\gamma(x)^{-1}\cdot \gamma (xz)= A_\g(x)^{-1}\cdot
n_{\gamma}^{-1}\cdot n_{\gamma}\cdot A_\g(xz)=A_\g(z)=
L_{n_\gamma^{-1}}\circ \gamma(z),\] where the third equality holds
as $A$ is an automorphism and $N/T=\tilde O$. This establishes the
above equality only on $\tilde O$, not on $\tilde F$, but both
sides of the equality make sense as elements of $H/I$, and since
$\tilde F$ is the fixed point set of $I$, any two elements of
$H/I$ that coincide on $\tilde O$ must coincide on $\tilde F$. Now
the right hand side is independent of $x$, which implies that the
rotational part of $\gamma$ is independent of $x$. This means that
the bundle $(\tilde O\times\mathbb R^k)/\pi_1(V)$ is a flat
$O(k)$-bundle, and hence so is the normal bundle of $F$ in $V$.
\end{proof}

\subsection{Normal bundle to a stratum is flat}
\label{subsec: norm bundle to strat} Let $F_t$ be the stratum
from Section~\ref{subsec: norm bundle in strat}.

\begin{lem} \label{lem: norm bundl to statum}
For all large $t$, the restriction to $O_t$ of the normal bundle
to $F_t$ in $H_t/\Gamma$ is flat.
\end{lem}
\begin{proof}
Let $\nu_t$ be the restriction to $O_t$ of the normal bundle to
$F_t$ in $H_t/\Gamma$. By an obvious contradiction argument it
suffices to show that any subsequence of $\{t_i\}$
has a subsequence for which $\nu_t$'s are flat; thus passing to
subsequences during the proof causes no loss of generality. Since
$k_t=\dim(\nu_t)$ can take only finitely many values, we pass to a
subsequence for which $k_t$ is constant; we then denote $k_t$ by
$k$. We look at the Grassmanian $G^k(TM)$ of $k$-planes in the
$TM$ with the metric induced by $g$. Of course, for $k$-planes
tangent to $H_t/\Gamma$ the metric is also induced by $g_t$. We
fix a point $o_t\in O_t$ and denote by $G^k_{t}$ the fiber of
$G^k(TM)$ over $o_t$. Let $l_t$ be the fiber of $\nu_t$ over
$o_t$.

The fibers of $G^k(TM)$ are (non-canonically) pairwise isometric
via the Levi-Civita parallel transport of $(M,g)$. Let $\rho\ll
diam(G^k_{t})$ be such that the center of mass of~\cite{Karcher}
is defined in any $4\rho$-ball of $G^k_{t}$. For large enough $t$,
the ball $B(l_t,4\rho)\subset G^k_{t}$ contains no $k$-plane
tangent to $F_t$ because $\rho\ll\diam(G^k_{t})$, while $l_t$,
$TF_t$ are $h_t$-orthogonal and $h_t$, $g_t$ are $\e$-close, so
that the distance between $l_t$ and any $k$-plane in $TF_t$ is
within $\e$ of $\diam(G^k_{t})$.

Denote by $\nabla^\infty$ the connection from Section~\ref{sec:
rot hom} associated with the parallel transport $P^\infty$. Also
denote by $\nabla_t$ the Levi-Civita connection of $h_t$, and let
$P_t$ be its parallel transport. The second fundamental form of
$H_t$ is uniformly bounded hence $P^\infty$ and $P_t$ are close
over any short loop in $H_t/\Gamma$.

Since $\nabla^\infty$ is flat and compatible with the metric $g$,
it defines the holonomy homomorphism
$\phi_t\co\G\to\Iso(T_{o_t}M)$. Let $R_t$ be the closure of
$\phi_t(\G)$ in $\Iso(T_{o_t}M)$. Since $R_t$ is compact, there
exists a {\it finite} subset $S_t\subset\Gamma$ so that for each
$r_t\in R_t$ there exists $s\in S_t$ such that $r_t(l_t)$,
$\phi_t(s)(l_t)$ are $\rho$-close in $G^k_{t}$. The direction
orthogonal to $H_t/\G$ is $\nabla^\infty$-parallel so the
$R_t$-action preserves the subspace $TH_t/\Gamma\subset TM$. Also
$P^\infty$ defines an isometry between $G^k_{t}$'s that is
$\phi_t$-equivariant and $R_t$-equivariant. Thus $S_t$ can be
chosen independently of $t$, and we denote $S_t$ by $S$.

If $t$ is sufficiently large, then for every $s\in S$, the
$k$-planes $\phi_t(s)(l_t)$ and $l_t$ are $\rho$-close in
$G^k_{t}$. Indeed, by the exponential convergence of geodesics, if
$t$ is large, then $s$ can be represented by a short loop based at
$o_t$. Since the inclusion $O_t\hookrightarrow H_t/\Gamma$ induces
a homotopy equivalence the loop can be assumed to lie on
$O_t\subset F_t$. Since $F_t$ is $h_t$-totally geodesic and $l_t$
is orthogonal to $F_t$ and has complementary dimension, $l_t$ is
fixed by $P_t$ along any loop in $F_t$ based at $o_t$. Hence the
same is almost true for $P^\infty$, provided the loop is short
enough, which proves the claim since $S$ is finite.

Thus, by the triangle inequality in $G^k_t$, the $k$-planes
$r_t(l_t)$ and $l_t$ are $2\rho$-close for all $r_t\in R_t$. Now
let $\bar l_t$ be the center of mass of all $r_t(l_t)$'s with
$r_t\in R_t$. Clearly, $\bar l_t$ is $2\rho$-close to $l_t$, hence
$\bar l_k$ is transverse to $TF_t$. Since $P^\infty$ preserves
$TH_t/\Gamma$, each $r_t(l_t)$, and hence $\bar l_t$, is tangent
to $H_t/\Gamma$. Now we translate $\bar l_t$ around $O_t$ using
$P^\infty$ along paths of length $\le\diam(O_t)<\e$. Since $\bar
l_t$ is $R_t$-invariant and $\nabla^\infty$ is flat, this gives a
well-defined $k$-dimensional flat $C^0$ subbundle $\bar\nu_t$ of
the restriction of $TH_t/\Gamma$ to $O_t\subset F_t$.

Note that $P_t$ takes any $k$-plane tangent to $F_t$ to a
$k$-plane tangent to $F_t$ since $F_t$ is totally geodesic. On the
other hand, $\bar l_t$ is uniform distance away from any $k$-plane
in $TF_t$, so its $P_t$-image remains far away from any $k$-plane
in $TF_t$. Since on short paths $P_t$ is close to $P^\infty$, we
conclude that $\bar\nu_t$ is transverse to $TF_t$ for large $t$.
Thus $\bar\nu_t$ is $C^0$ isomorphic to $\nu_t$, so that $\nu_t$
carries a $C^0$ flat connection.

In fact, $\nu_t$ also carries a smooth flat connection. Indeed, in
the universal cover $\tilde O_t$, the $C^0$ flat connection
defines a $\G$-equivariant $C^0$-isomorphism of the pullback of
$\nu_t$ to $\tilde O_t$ onto $\tilde O_t\times\Rk$, where $\G$
acts as the covering group on the $\tilde O_t$-factor and via a
holonomy homomorphism on the $\Rk$-factor. This is a smooth
action, so the quotient $(\tilde O_t\times\Rk)/\G$ is a smooth
flat vector bundle that is $C^0$ isomorphic to $\nu_t$. But any
$C^0$-isomorphic bundles are smoothly isomorphic because the
continuous homotopy of classifying maps can be approximated by a
smooth homotopy.
\end{proof}

\section{Infranilmanifolds are horosphere quotients}

Z.~Shen constructed in~\cite{Shen}
a pinched negatively curved warped product
metric on the product of
an arbitrary infranilmanifold and $(0,\infty)$ so that
the metric is complete
near the $\infty$-end, but is incomplete at the $0$-end. Here we
modify Shen's construction to produce a complete pinched
negatively curved metric on the product of any infranilmanifold
with $\mathbb R$.

Let $G$ be a simply-connected
nilpotent Lie group acting on itself by left
translations, and let $K$ be a compact subgroup of
$\mathrm{Aut}(G)$, so that the semidirect product $G\rtimes K$
acts on $G$ by affine transformations. Taking product with
the trivial $G\rtimes K$-action on $\mathbb R$, we get
a $G\rtimes K$-action on $G\times\mathbb R$ for which we
prove the following.

\begin{thm}\label{thm: warp infranil}
$G\times\mathbb R$ admits a complete $G\rtimes K$-invariant
Riemannian metric of pinched negative curvature. In particular, if
$N$ is an infranilmanifold, then $N\times\mathbb R$ carries a
complete metric of pinched negative curvature.
\end{thm}
\begin{proof}
The Lie algebra $L(G)$ can
be written as
\[L(G)=L_0\supset L_1\supset\cdots\supset L_k\supset L_{k+1}=0\]
where $L_{i+1}=[L_0,L_i]$.
Note that $[L_i,L_j]\subset L_{i+j+1}$. Indeed, assume
$i\le j$ and argue by induction on $i$. The case $i=0$ is obvious
and the induction step follows from the Jacobi identity
and the induction hypothesis, because
$[L_i,L_j]=[[L_0,L_{i-1}],L_j]$ lies in
\[\mathrm{span}([[L_{i-1},L_j],L_0], [[L_0,L_j],L_{i-1}])\subset
\mathrm{span}([L_{i+j},L_0], [L_{j+1},L_{i-1}])= L_{i+j+1} \]
The group $K$ preserves each $L_i$, so we can choose a
$K$-invariant inner product $\langle\ ,\ \rangle_0$ on $L$. Let
\[F_i=\{X\in L_i\co \langle X,Y\rangle_0=0\ \mathrm{for}\ Y\in L_{i+i}\}.\]
Then $L=F_0\oplus\cdots\oplus F_k$. Define a new $K$-invariant inner
product $\langle\ ,\ \rangle_r$ on $L$ by $\langle
X,Y\rangle_r=h_i(r)^2\langle X,Y\rangle_0$ for $X,Y\in F_i$, and
$\langle X,Y\rangle_r=0$ if $X\in F_i$, $Y\in F_j$ for $i\neq j$,
where $h_i$ are some positive function defined below. This defines
a $G\rtimes K$-invariant Riemannian metric $g_r$ on $G$.

Let $\alpha_i=i+1$ with $i=0,\cdots, k$ and $a=k+1$. Now define the
warping function $h_i$ to be a positive, smooth, strictly
convex, decreasing
function that is equal to $e^{-\a_i r}$ if $r\ge 1$, and is equal to
$e^{-ar}$ if $r\le -1$; such a function exists since
$a\ge a_i$ for each $i$.
Thus $h_i^\prime<0<h_i^{\prime\prime}$, and the functions
$\frac{h_i^\prime}{h_i}$,
$\frac{h_i^{\prime\prime}}{h_i}$ are uniformly bounded
away from $0$ and $\infty$.

Define the warped product metric on $G\times\mathbb R$ by $g=s^2
g_r+dr^2$, where $s>0$ is a constant; clearly $g$ is a complete
$G\rtimes K$-invariant metric.
A straightforward tedious computation
(mostly done e.g. in~\cite{BW}) yields for $g$-orthonormal vector
fields $Y_s\in F_s$ that
\begin{eqnarray*}
 & \langle R_g(Y_i,Y_j) Y_j,Y_i\rangle_g=
\frac{1}{s^2}\langle R_{g_r}(Y_i,Y_j) Y_j,Y_i\rangle_{g_r}-
\frac{h_i^\prime h_j^\prime}{h_ih_j},\\
& \langle R_g(Y_i,Y_j) Y_l,Y_m\rangle_g= \frac{1}{s^2}\langle
R_{g_r}(Y_i,Y_j) Y_l,Y_m\rangle_{g_r}
\ \ \ \mathrm{if}\ \{i,j\}\neq \{l,m\},\\
 & \langle R_g(Y_i,\frac{\partial}{\partial r})
\frac{\partial}{\partial r}), Y_i \rangle_g=
-\frac{h_i^{\prime\prime}}{h_i},\ \ \ \ \ \langle
R_g(Y_i,\frac{\partial}{\partial r})
\frac{\partial}{\partial r}), Y_j \rangle_g=0\ \ \ \mathrm{if}\ i\neq j, \\
& \langle R_g(\frac{\partial}{\partial r},Y_i) Y_j,Y_l\rangle_g=
\left(\frac{h_j^\prime}{2h_j}+\frac{h_l^\prime}{2h_l}\right)
\left(\langle [Y_j,Y_i],Y_l\rangle_g+\langle
[Y_i,Y_l],Y_j\rangle_g +\langle [Y_j,Y_l],Y_i\rangle_g\right).
\end{eqnarray*}

{\bf Correction} (added on August 28, 2010): {\it
The above formula
for $\langle R_g(\frac{\partial}{\partial r},Y_i) Y_j,Y_l\rangle_g$
is incorrect. A correction can be found in
Appendix C of~\cite{Bel-rh-warp} where it is explained
why the mistake does not affect other results
of the present paper.}

Since $[L_i,L_j]\subset L_{i+j+1}$, we have
for $Z=\sum_{i=0}^k Z_i$ and $W=\sum_{j=0}^k W_j$ with $Z_i, W_i\in F_i$
\[
|[Z,W]|_{g_r}\le \sum_{ij}|[Z_i,W_j]|_{g_r}\le\sum_{ij}\sum_{s>i+j} h_s |[Z_i,W_j]|_{g_0}
\]
The above choice of $a_i$'s implies that if $r\ge 1$,
then $\sum_{s>i+j} h_s\le k h_i h_j$. Also
$|[Z_i,W_j]|_{g_0}\le C |Z_i|_{g_0}|W_j|_{g_0}$ where $C$ only depends on
the structure constants of $L$, so that we conclude
\[
|[Z,W]|_{g_r}\le Ck |Z_i|_{g_0}|W_j|_{g_0} \sum_{ij} h_ih_j\le Ck(k+1)|Z|_{g_r}|W|_{g_r}.
\]
It follows that if $r\ge 1$, then the norm of the curvature tensor
of $g_r$ is bounded in terms of $C$, $k$~\cite[Proposition 3.18]{Cheeger-Ebin}.
The same conclusions
trivially hold for $r\le -1$, because then $g_r$ is the rescaling
of $g_0$ by a constant $e^{-ar}>1$, and also for $r\in [-1,1]$ by
compactness, since $g_r$ is left-invariant and depends continuously
of $r$. Hence $\langle R_g(Y_i,Y_j) Y_l,Y_m\rangle_g\to 0$ as $s\to \infty$.

Also $\langle R_g(\frac{\partial}{\partial r},Y_i) Y_j,Y_l\rangle_g\to 0$
as $s\to\infty$, because
\[|\langle [Y_j,Y_i],Y_l\rangle_g|=
s^2|\langle [Y_j,Y_i],Y_l\rangle_{g_r}|\le
s^2C(k+1)|Y_j|_{g_r}|Y_i|_{g_r}|Y_l|_{g_r}\le C(k+1)/s,\]
where the last inequality holds since
$s|Y|_{g_r}=1$ for any $g$-unit vector $Y$.
It follows that as $s\to\infty$, then $R_g$ uniformly converges to
a tensor $\bar R$ whose nonzero components are
\[ \bar R(Y_i,Y_j, Y_j,Y_i)=
-\frac{h_i^\prime h_j^\prime}{h_ih_j}\ \ \mathrm{and}\ \ \bar
R\left(Y_i,\frac{\partial}{\partial r}, \frac{\partial}{\partial r},
Y_i\right)= -\frac{h_i^{\prime\prime}}{h_i}.\]
Thus $g$ has pinched negative curvature for all large $s$.
\end{proof}

\begin{cor}
Let $E$ be the total space of a flat  Euclidean vector bundle over an
infranilmanifold $I$. Then $E$ is infranil, in particular,
$E\times \mathbb R$ admits a complete
Riemannian metric of pinched negative curvature.
\end{cor}
\begin{proof}
Fix a flat  Euclidean  $\mathbb R^k$-bundle over the infranilmanifold $I$, and write
$I$ as $G_0/\Gamma$ where $G_0$ a simply-connected nilpotent Lie
group and $\Gamma$ is a discrete cocompact group of affine
transformation of $G_0$ that acts freely. Look at the nilpotent
group $G=G_0\times\mathbb R^k$, and let $\Gamma$ act on the
$\mathbb R^k$-factor via the holonomy of the flat bundle
$\Gamma\cong\pi_1(I)\to O(k)$. Then the infranilmanifold
$G/\Gamma$ is diffeomorphic to the total space of the flat bundle
we started with. By Theorem~\ref{thm: warp infranil},
$G/\Gamma\times \mathbb R$ carries a complete metric of pinched
negative curvature.
\end{proof}

\section{On geometrically finite manifolds}\label{section: manifolds}

\begin{proof}[Proof of Corollary~\ref{cor: geom finite}]
Let $X/\G$ be a geometrically finite pinched negatively
curved manifold, let $\Omega$ be the domain of discontinuity
and $\L$ be the limit set for the $\G$-action at infinity.
Let $C_\e$ be the $\e$-neighborhood of the convex hull of
$L$. Then $C_\e/\G$ is a codimension zero $C^1$ submanifold
of $X/\G$ that is homeomorphic to $(X\cup\Omega)/\G$ by
pushing along geodesic rays orthogonal to $\partial C_\e/\G$.
This homeomorphism restricts to diffeomorphism
on the interiors $X/\G\to \mathrm{Int}(C_e)/\G$.

By the discussion in~\cite[pp263-264]{Bow2}, for each end of
$(X\cup\Omega)/\G$ there is a parabolic subgroup $\G_z\le \G$
stabilizing a point $z\in\partial_\infty X$ such that the
end has a neighborhood homeomorphic  to a neighborhood of the
unique end of $(X\cup\partial_\infty X\setminus\{z\})/\G_z$. Again,
this homeomorphism restricts to diffeomorphism on the interiors.

By pushing along trajectories of Busemann flow,
$(X\cup\partial_\infty X\setminus\{z\})/\G_z$
is homeomorphic to the $\G_z$-quotient of a closed
horoball $H_z$ centered at $z$. Note that $H_z/\G_z$
is a $C^2$ submanifold of $X/\G_z$, and $H_z/\G_z$
is $C^1$ diffeomorphic to the product of $[0,\infty)$ and
a horosphere quotient, that
by Theorem~\ref{main thm} is diffeomorphic
to the interior of a compact manifold $L_z$. So $H_z/\G_z$
is diffeomorphic to the
interior of $L_z\times [0,1]$, in which we smooth corners.
Compactifying each $z$-end of $C_\e/\G$
with $L_z\times [0,1]$, we get a compact $C^1$ manifold
whose interior is diffeomorphic to $X/\G$.
\end{proof}

\appendix
\section{Lemmas on nilpotent groups}
\label{app: nilp}
\begin{lem}\label{root}
Given a finitely generated nilpotent group $\G$ and a positive
integer $n$, there exists a finite subset $F\subset \G$ such that for
any $g\in \G$ there is $f\in F$ and $x\in \G$ with $gf=x^n$.
\end{lem}
\begin{proof}
We argue by induction on the nilpotency degree of $\G$. If $\G$ is
abelian, then the $n$th powers of elements of $\G$ form a finite
index subgroup, and we can take $F$ to the set of coset
representatives of this subgroup. In general, if $Z$ denotes the
center of $\G$, then by induction the result is true for $\G/Z$
for the finite subset $\{aZ: a\in F_1\}$ of $G/Z$, where $F_1$ is
some finite subset of $\G$. Thus, an arbitrary $g\in \G$ satisfies
$gf_1=x^n z$ for some  $f_1\in F_1$, $x\in \G$, $z\in Z$. Again
since $Z$ is abelian, the set $Z^n$ of $n$th powers is a finite
index subgroup of $Z$. Let $F_2$ be a set of coset representatives
of $Z^n$ in $Z$ so that $z=y^n f_2$ for some $y\in Z$, $f_2\in
F_2$. Then $gf_1=x^ny^nf_2=(xy)^nf_2$, so $g f_1 f_2^{-1}
=(xy)^n$, and the assertion holds for $\G$ with $F=F_1F_2^{-1}$.
\end{proof}

\begin{lem}\label{root2}
Let $\G$ be a finitely generated nilpotent group. Then there
exists a finite set $F\subset \G$ such that for any
$x\in \G$, $g\in [\G,\G]$ there is $h\in [\G,\G]$, $f\in F$
with $x^{2}g=(x h)^{2}f$.\end{lem}
\begin{proof}
As usual we denote $\G_1=\G$, $\G_{i+1}=[\G,\G_i]$, so that the
nilpotency degree $k$ of $\G$ is the largest integer for which
$\G_k$ is nontrivial. Since $[\G,\G_{k}]$ is trivial, $\G_{k}$
lies in the center of $\G$. We argue by induction on $k$. The case
$k=1$, i.e. when $\G$ is abelian, is obvious for $F=\{1\}$. If
$\G$ is of nilpotency degree $k>1$, then by induction the
statement is true in $\G/\G_{k}$. Let $F_1\subset \G$ be the set
of coset representatives of the corresponding finite set for
$\G/\G_{k}$, so that given $x\in\G$, there exists $z\in\G_{k}$,
$g,h\in [\G,\G]$, $f_1\in F_1$ with $x^{2}g=(x h)^{2}f_1z$. By
Lemma~\ref{root} applied to $\G_k$, we get $z=y^{2}f_2$ for some
$y\in \G_k$, $f_2\in F_2$ where $F_2$ is a finite subset of
$\G_k$. Then $x^{2}g=(x h)^{2}f_1z=(x y h)^{2}f_1f_2$, and since
$yh\in [\G,\G]$ the proof is complete.
\end{proof}

\begin{lem}
\label{lem: max comp}
Let $C$ be a maximal compact subgroup of a
connected nilpotent Lie group $N$. Then
$C$ is equal to the unique maximal compact subgroup
of the center of $N$, in particular $C$ is a torus.
\end{lem}
\begin{proof}
Any maximal compact subgroup $C$ of $N$ lies in the center $Z$ of
$N$~\cite{Iwa}.
Hence $C$ is also maximal compact in $Z$. A maximal
compact subgroup is homotopy equivalent to the ambient group, hence
since $N$ is connected, so is $Z$ and $C$. Now $Z$ connected
abelian, hence $Z$ is isomorphic to the product of a real vector space
and a torus, so $C$ equals to the torus.
\end{proof}

\section{Isometries are smooth}
\label{apendix: iso are c3}

\begin{prop}
Let $X$ be a smooth manifold equipped with a complete $C^0$-Riemannian
metric of curvature bounded above and below in the comparison sense.
Then the isometry group acts on $X$ by $C^3$-diffeomorphisms.
\end{prop}
\begin{proof}
The isometry group of any complete locally-compact metric space is
locally compact~\cite{KN}.
By~\cite[Chapter 5]{MZ}, any locally compact subgroup
of $\mathrm{Diffeo}^r(X)$ with $r>0$ is a Lie group and the action is $C^r$.
Thus, it suffices to show that each individual isometry is $C^3$.
The construction of harmonic coordinates in~\cite{Nik2} starts
with the $C^0$-distance coordinates $(d(x,a_1),...,d(x,a_n))$
at $x\in X$, where the geodesic segments $[x, a_i]$ are pairwise
orthogonal, and then solves the Dirichlet problem in a small ball
around $x$ with values on the boundary sphere given by $d(x,a_i)$.
The solutions are the so-called harmonic coordinates. Their transition
functions are $C^{3,\a}$ (and the metric tensor in this coordinates
is $C^{1,\a}$ even though we do not need this fact here).
This construction is clearly invariant under isometries, so
any isometry has the same smoothness in harmonic coordinates
as the identity map, namely $C^{3,\a}$.
\end{proof}

\section{Local formula of Ballmann and Br\"{u}ning}\label{BBB}

Let $X$ be a simply-connected manifold of pinched negative
curvature. Fix a point at infinity of $X$, and let $T$ be the unit
vector field tangent to the Busemann flow $b_t(x)$ towards that point.
For a curve
$\a(s)$ in $X$, denote $\a(0)=x$ and $\a'(0)=u$. Look at the
$1$-parameter family of geodesic rays $\a(s,t)=b_t(\a(s))$
and the corresponding family of Jacobi fields $J$.
Let $v,w\in T_x X$ and
let $X(t,s)$, $Y(t,s)$ be vector fields along $a(s,t)$ such that
$X(0,0)=v$, $Y(0,0)=w$ and $\nabla_T X=\nabla_T Y=0$. Then one defines
a tensor field $\bar S$  by
\[
\langle \bar{S}(u,v),w\rangle=-\int_0^\infty \langle
R(T,J)X,Y\rangle(t,0) dt
\]
Ballmann and Br\"uning~\cite{BB} define a new connection
$\bar{\nabla}$ by $\bar{\nabla}=\nabla-\bar{S}$, and show that
$\bar{\nabla}$ is a $C^0$ flat connection that is compatible with
the metric, and that satisfies $\bar{\nabla}T=0$ and
\begin{equation}
\label{xy}
|\bar{\nabla}_XY-\nabla_XY|\le C(a)|X||Y| \text{ for any }X,Y
\end{equation}

Since $\bar{\nabla}$ is flat, for any $x,y\in M$ we
have a well defined parallel transport with respect to
$\bar{\nabla}$ from $x$ to $y$ which we denote by
$P^{\bar{\nabla}}_{xy}$.

\begin{lem}\label{nabla=infty}
For any $x,y\in M$ parallel transport through infinity
$P^\infty_{xy}$ coincides with $P^{\bar{\nabla}}_{xy}$.
\end{lem}
\begin{proof}
First suppose $x$, $y$ lie in the same horosphere $H_{t_0}$.
Consider the quadrangle $x b_t(x)b_t(y)y$. By flatness
$P^{\bar{\nabla}}_{xy}$ is equal to the $\bar{\nabla}$
parallel transport along three other sides of this quadrangle
$P^{\bar{\nabla}}_{b_t(y)y}\circ
P^{\bar{\nabla}}_{b_t(x)b_t(y)}\circ
P^{\bar{\nabla}}_{xb_t(x)}$. Also if $\a$ is a trajectory of the
Busemann flow, then $J=T$ so that $\bar S=0$,
therefore parallel transports along
Busemann trajectories coincide with Levi-Civita parallel
transports. Thus
\[
P^{\bar{\nabla}}_{xy}=P_{b_t(y)y}\circ
P^{\bar{\nabla}}_{b_t(x)b_t(y)}\circ P_{xb_t(x)}
\]

Since $d(b_t(x)b_t(y))\to 0$ as $t\to\infty$, by~(\ref{xy}) we
have that $P^{\bar{\nabla}}_{b_t(x)b_t(y)}$ becomes arbitrary
close to $P_{b_t(x)b_t(y)}$ for large $t$ and therefore
\begin{equation}\label{nablap}
P^{\bar{\nabla}}_{xy}=\lim_{t\to\infty}P_{b_t(y)y}\circ
P_{b_t(x)b_t(y)}\circ P_{xb_t(x)}
\end{equation}

Similarly, by construction $P^\infty$ commutes with Busemann flow
and hence
\[
P^{\infty}_{xy}=P_{b_t(y)y}\circ P^{\infty}_{b_t(x)b_t(y)}\circ
P_{xb_t(x)}
\]
As before $P^{\infty}_{b_t(x)b_t(y)}$ becomes
arbitrary close to $P_{b_t(x)b_t(y)}$ for large $t$ and therefore
\begin{equation}\label{inftyp}
P^{\bar{\nabla}}_{xy}=\lim_{t\to\infty}P_{b_t(y)y}\circ
P_{b_t(x)b_t(y)}\circ P_{xb_t(x)}
\end{equation}
Comparing~(\ref{nablap}) and~(\ref{inftyp}) we conclude that
$P^{\infty}_{xy}=P^{\bar{\nabla}}_{xy}$ for any $x$, $y$
in the same horosphere.

Finally, since $\bar\nabla$ and $\nabla^\infty$ are flat, and
their parallel transports coincide with the Levi-Civita parallel
transport along trajectories of the Busemann flow
$P^{\infty}_{xy}=P^{\bar{\nabla}}_{xy}$ for any $x$, $y$.
\end{proof}

\section{Concave functions  and submetries on Alexandrov spaces}
\label{A:crit}
The proof of the following lemma is due to A.~Petrunin.
\begin{lem}\label{crit}
Let $X$ be an Alexandrov space of $\curv\ge k$ with $\partial X=\emptyset$. Let $f\co X\to
\mathbb R$ be a  Lipschitz function  with a local maximum at $q$.
Suppose $f$ is strictly concave on an open set $U$ containing $q$.
Then $d(\cdot, q)$ has no critical points on $U\backslash \{q\}$.
\end{lem}
\begin{proof}
Let $x\in U, x\ne q$. By~\cite{PP} $\nabla f(x)$ is
defined to be equal to $v\in T_xX$ if $df(v)=|v|^2$ and
$\frac{df(u)}{|u|}$ attains a  positive maximum at $v$. Since $x$
is not a point of maximum of $f$, $\nabla f(x)\ne 0$ and $|\nabla
f(z)|\ge c>0$ for  all $z$ near $x$.
Consider the gradient flow  for $f$
as defined in~\cite{PP}. Consider  a gradient line $\g(t)$ passing
through $x$ so that $\g(0)=x$. By~\cite{PP} the curve $\g(t)$ is
locally Lipschitz.

We claim that $\g(t)$ can be extended to be a gradient line  
defined on $(-\e,\infty)$ for some $\e>0$.
Indeed, suppose this is not true. 
By ~\cite[Lemma 3.2.1(a)]{PP} the gradient flow of {\it any}
concave function is $1$-Lipschitz, so it defines 
a deformation retraction of any superlevel set
$\{f\ge c\}$ contained in $U$ onto $q$ 
and hence  $\{f\ge c\}$ is contractible 
(for a different proof see also~\cite[Lemma 5.2]{Kap}).  
Take $c=f(x)$ and let $\e>0$ be small enough 
so that $\{f\ge c-\e\}$ is contained in $U$. 
Since $\g$ can not be extended backwards beyond zero,  
the gradient flow gives a deformation retraction of 
$\{f\ge c-\e\}\backslash \{x\}$ onto $q$. 
To see that this is impossible we prove that 
$H^{n-1}(\{f\ge c-\e\}\backslash \{x\},\Z_2)\cong \Z_2$.

Indeed, since $\{f\ge c-\e\}$ is contractible, from the long  exact cohomology sequence of a pair we see that $$H^{n-1}(\{f\ge c-\e\}\backslash \{x\},\Z_2)\cong H^n(\{f\ge c-\e\},\{f\ge c-\e\}\backslash \{x\},\Z_2).$$

By ~\cite{Per2}, for some very small $\d>0$, the ball
$B(x,\d)$ is contractible and is contained in
$\{f\ge c-\e\}$, and $d(\cdot, x)$ has no critical points in 
$B(x,2\d)\backslash\{x\}$. Therefore, by excision
\[
H^n(\{f\ge c-\e\},\{f\ge c-\e\}\backslash \{x\},\Z_2)\cong 
H^n(B(x,\d), B(x,\d)\backslash \{x\}, \Z_2)\]
which is isomorphic to 
$H^{n-1}(S(x,\d),\Z_2)\cong \Z_2$,
where the last equality holds since by~\cite{Per2},  
$S(x,\d)$ is  homotopy equivalent to the Alexandrov
space $\Sigma_{x}X$ with
$\partial \Sigma_{x}X=\emptyset$, and since any Alexandov space
without boundary has top-dimensional $\Z_2$-cohomology isomorphic
to $\Z_2$~\cite{GrPe}.

Let $v$ be any {\it left}  tangent of $\g$ at $0$. Here
following~\cite{PP} we say that $v\in T_xX$ is a left tangent
vector if $v=\lim_{i\to\infty}\frac{1}{|t_i|}exp_x^{-1}\g(t_i)$
for some sequence $t_i\to 0^-$.
(One should think of $v$ as $-\nabla f(x)$). By above
$f(\g(t))-ct$ is nondecreasing for small $t$. Therefore $v\ne 0$.
We claim that $\angle uv\ge\pi/2$ for any direction $u$ of a
shortest geodesic from $x$ to $y$ with $f(y)> f(x)$. Indeed,
by~\cite[Lemma 3.2.1(a)]{PP} the gradient flow of {\it any}
concave function is $1$-Lipschitz. Consider a cutoff concave
function $\hat{f}(\cdot)=\min\{f(\cdot), f(y)\}$. Clearly the
forward gradient flow of $\hat{f}$ fixes $y$ and coincides with
the gradient flow of $f$ near $x$. Since the gradient flow of
$\hat{f}$ is $1$-Lipschitz,  $d(y,\g(t))$ is nonincreasing near
$t=0$. By the first variation formula, this implies that $\angle
uv\ge\pi/2$. In particular, this is true for $y=q$. Since any
shortest from $x$ to $q$ points strictly inside the convex set
$\{f\ge f(x)\}$  this inequality is in fact strict, i.e $\angle
uv>\pi/2$.
\end{proof}

\begin{defn}
A map $f\co X\to Y$ between two metric spaces is called a submetry
if for any $x\in X$ and any $r>0$ one has $f(B_r(x))=B_r(f(x))$.
\end{defn}
We need some properties of submetries collected below.
\begin{lem}\label{sub}
Suppose $X$ is an Alexandrov space of $curv\ge k$ and $f\co X\to
Y$ is a submetry. Then
\begin{enumerate}[a)]
\item $Y$ is an Alexandrov space of $curv\ge k$~\cite{BGP} .
\item
For any $x\in X, y\in Y$ we have $d(x,f^{-1}(y))=d(f(x),y)$.
\item
For any $x\in X$ and any shortest geodesic $\g\co [0, 1]\to Y$
with $\g(0)=f(x)$ there exists a shortest geodesic 
$\tilde{\g}\co [0, 1]\to Y$ with $\tilde{\g}(0)=x$ such that
$f(\tilde{\g}(t))=\g(t)$. The geodesic $\tilde{\g}$ is called a
horizontal lift of $\g$. Moreover, if $\g(t)$ can be extended to a
shortest $\g\co [-\e, 1]\to Y$ then the horizontal lift is unique.
\end{enumerate}
Moreover the statement a) is true locally in the following sense:
if $B(p,1)\subset U\subset X$ and $g\co U\to Y$ is a submetry and
$\diam (g^{-1}(y))<1/10$ for any $y\in Y$, 
then for any triangle with vertices in
$B(g(p),1/4)$ the Toponogov comparison with $curv\ge k$ holds.

\end{lem}
The proofs of parts b) and c) are elementary and
left to the reader (see~\cite{Lyt} for details).

\bibliographystyle{amsalpha}
\bibliography{nilp-final-withcomments}
\end{document}